\catcode\endlinechar = 10%
%
\voffset = 4.6mm%
\vsize = 237mm%
\hoffset = 4.6mm%
\hsize = 150mm%
\font\fourteenrm = cmr12 at 14pt%
\font\twelverm = cmr12%
\font\eightpointfourrm = cmr8 at 8.4pt%
\font\sixrm = cmr6%
\font\twelvebf = cmbx12%
\font\eightpointfourbf = cmbx8 at 8.4pt%
\font\sixbf = cmbx6%
\font\twelveit = cmti12%
\font\eightpointfourit = cmti8 at 8.4pt%
\font\sixit = cmti7 at 6pt%
\font\twelvei = cmmi12%
\font\eightpointfouri = cmmi8 at 8.4pt%
\font\sixi = cmmi6%
\font\twelvesy = cmsy10 at 12pt%
\font\eightpointfoursy = cmsy8 at 8.4pt%
\font\sixsy = cmsy6%
\font\twelveex = cmex10 at 12pt%
\font\twentysc = cmcsc10 at 20pt%
\font\fourteensc = cmcsc10 at 14pt%
%
%
\font\twelvett = cmtt12 at 12pt%
\font\eightpointfourtt = cmtt8 at 8.4pt%
\font\sixtt = cmtt8 at 6pt%
\textfont0 = \twelverm%
\scriptfont0 = \eightpointfourrm%
\scriptscriptfont0 = \sixrm%
\textfont1 = \twelvei%
\scriptfont1 = \eightpointfouri%
\scriptscriptfont1 = \sixi%
\textfont2 = \twelvesy%
\scriptfont2 = \eightpointfoursy%
\scriptscriptfont2 = \sixsy%
\textfont3 = \twelveex%
\scriptfont3 = \twelveex%
\scriptscriptfont3 = \twelveex%
\newfam\itfam%
\textfont\itfam = \twelveit%
\scriptfont\itfam = \eightpointfourit%
\scriptscriptfont\itfam = \sixit%
\newfam\bffam%
\textfont\bffam = \twelvebf%
\scriptfont\bffam = \eightpointfourbf%
\scriptscriptfont\bffam = \sixbf%
\newfam\ttfam%
\textfont\ttfam = \twelvett%
\scriptfont\ttfam = \eightpointfourtt%
\scriptscriptfont\ttfam = \sixtt%
\def\rm{\fam=0\twelverm}%
\def\it{\fam=\itfam\twelveit}%
\def\bf{\fam=\bffam\twelvebf}%
\baselineskip = 14.4pt%
\normalbaselineskip = \baselineskip%
\def\makefootline{%
	\baselineskip 28.8pt%
	\lineskiplimit 0pt%
	\line{\the\footline}%
}%
\rm%
\parindent = 1cm%
\parskip = 5pt%
\def\monthname{%
	\ifcase\month%
		\or January%
		\or February%
		\or March%
		\or April%
		\or May%
		\or June%
		\or July%
		\or August%
		\or September%
		\or October%
		\or November%
		\or December%
	\fi%
}%
\footline = {\hss\rm\folio\hss}%
\def\Document{article}%
\def\Title#1{{%
	\parskip = 0pt%
	\baselineskip = 24pt%
	{%
		\twentysc%
		\halign{%
			\hbox to \hsize{\hfill##\hfill}\cr%
			#1\cr%
		}%
	}%
	\centerline{\fourteenrm B.R.S.~Roso}%
	\baselineskip = 16.8pt%
	\centerline{\fourteenrm\monthname\ \number\year}%
	\centerline{\fourteenrm University College London}%
	\vskip 20pt%
}}%
\newcount\SectNum%
\newcount\EltNum%
\SectNum = 0%
\def\Abs{{%
	\vskip\baselineskip%
	\fourteensc%
	\centerline{Abstract}%
}}%
\def\Sect#1{%
	\advance\SectNum by 1%
	\EltNum = 1%
	{%
		\vskip\baselineskip%
		\goodbreak%
		\fourteensc%
		\centerline{\number\SectNum.\enskip#1}%
		\nobreak%
	}%
}%
\long\def\Elt#1#2{%
	{%
		\goodbreak%
		\par%
		\noindent%
		\def\EltName{%
			\number\SectNum%
			.%
			\number\EltNum%
		}%
		\def\EltFullName{%
			#1%
			\ 
			\EltName
		}%
		{\bf \EltFullName:\ }%
		{\ignorespaces#2}%
		\goodbreak%
		\par%
	}%
	\advance\EltNum by 1%
}%
\def\Label#1{%
	\global\edef#1{{\EltFullName}}%
}%
\long\def\Thm#1{%
	\Elt{Theorem}{#1}%
}%
\long\def\Prop#1{%
	\Elt{Proposition}{#1}%
}%
\long\def\Lem#1{%
	\Elt{Lemma}{#1}%
}%
\long\def\Def#1{%
	\Elt{Definition}{#1}%
}%
\long\def\Rem#1{%
	\Elt{Remark}{#1}%
}%
\long\def\Eg#1{%
	\Elt{Example}{#1}%
}%
\long\def\Cor#1{%
	\Elt{Corollary}{#1}%
}%
\long\def\Conj#1{%
	\Elt{Conjecture}{#1}%
}%
\long\def\Quest#1{%
	\Elt{Question}{#1}%
}%
\def\qed{{%
	\unskip%
	\nobreak%
	\hfil%
	\penalty50%
	\hskip2em%
	\hbox{}%
	\nobreak%
	\hfil%
	{\bf QED}%
	\parfillskip=0pt%
	\finalhyphendemerits=0%
	\par%
}}%
\long\def\Proof#1{%
	{%
		\goodbreak%
		\par%
		\noindent%
		{\bf Proof:\ }%
		{\ignorespaces #1}%
		\qed%
	}%
	\goodbreak%
}%
\long\def\Ack#1{%
	{%
		\goodbreak%
		\par%
		\noindent%
		{\bf Acknowledgements:\ }%
		{\ignorespaces #1}%
	}%
	\goodbreak%
}%
\long\def\Keywords#1{%
	{%
		\goodbreak%
		\par%
		\noindent%
		{\bf Keywords:\ }%
		{\ignorespaces #1}%
	}%
	\goodbreak%
}%
\long\def\MSC#1{%
	{%
		\goodbreak%
		\par%
		\noindent%
		{\bf MSC2020:\ }%
		{\ignorespaces #1}%
	}%
	\goodbreak%
}%
%
%
\def\Fig#1{%
	{%
		\goodbreak%
		\par%
		\vskip\parskip%
		\def\EltName{%
			\number\SectNum.\number\EltNum%
		}%
		\def\EltFullName{%
			Figure\ \EltName%
		}%
		\hbox to \hsize {%
			\hfill%
			\vbox{%
				\halign{%
					\hfill%
					##%
					\hfill%
					\nobreak%
					\cr%
					#1%
					\cr%
					{\bf \EltFullName}%
					\cr%
				}%
			}%
			\hfill%
		}%
		\par%
	}%
	\advance\EltNum by 1%
}%

%
\font\twelvefrak = eufm10 at 12pt%
\font\eightpointfourfrak = eufm8 at 8.4pt%
\font\sixfrak = eufm6 at 6 pt%
\newfam\frakfam%
\textfont\frakfam = \twelvefrak%
\scriptfont\frakfam = \eightpointfourfrak%
\scriptscriptfont\frakfam = \sixfrak%
\def\frak{\fam\frakfam\twelvefrak}%
\font\twelverunes = frucms at 12pt%
\font\eightpointfourrunes = frucms at 8.4pt%
\font\sixrunes = frucms at 6pt%
\newfam\runesfam%
\textfont\runesfam = \twelverunes%
\scriptfont\runesfam = \eightpointfourrunes%
\scriptscriptfont\runesfam = \sixrunes%
%
%
%
\def\C{{\bf C}}%
\def\R{{\bf R}}%
\def\N{{\bf N}}%
\def\Q{{\bf Q}}%
\def\Z{{\bf Z}}%
%
%
\def\defop#1#2{\def#1{\mathop{#2{}}\nolimits}}%
\def\d{{\rm d}}%
\defop\D{{\rm D}}%
\def\conj#1{{\overline{#1}}}%
\defop\eval{{\rm ev}}%
\def\Hodge{\mathord{\ast}}%
\defop\ahd{\conj{\partial}}%
\defop\ad{{\rm ad}}%
\defop\Ad{{\rm Ad}}%
\def\Dirac{{\frak D}}%
\def\morphnamed#1{{%
	\setbox0 = \hbox{$\scriptstyle #1$}%
	\dimen0 = \wd0%
	\advance\dimen0 by 10pt%
	\mathop{%
		\;%
		\hbox to \dimen0 {\rightarrowfill}%
		\;%
	}\limits^{\box0}%
}}%
\def\converge#1{{%
	\setbox0 = \hbox{$\scriptstyle #1$}%
	\dimen0 = \wd0%
	\advance\dimen0 by 10pt%
	\mathop{%
		\;%
		\hbox to \dimen0 {\rightarrowfill}%
		\;%
	}\limits_{\box0}%
}}%
\defop\sign{{\rm sign}}%
\defop\index{{\rm ind}}%
\defop\And{\enskip\&\enskip}%
\defop\PD{{\rm PD}}%
\defop\LieDer{{\cal L}}%
\defop\cliff{{c\ell}}%
%
%

%
%
%
%
\def\U{{\rm U}}%
%
%
%
\defop\Hom{{\rm Hom}}%
\defop\Ker{{\rm Ker}}%
\defop\Image{{\rm Im}}%
\defop\Graph{{\rm Graph}}%
\defop\Coker{{\rm Coker}}%
\defop\Endo{{\rm End}}%
\defop\Sym{{\rm Sym}}%
\def\Homology{{\rm H}{}}%
\def\RedHomology{\tilde{\rm H}{}}%
\def\GeomCompl{{\rm c}}%
\def\BorelHomology{\GeomCompl\RedHomology}%
\def\HMFrom{\widehat{\rm HM}{}}%
\def\HMTilde{\widetilde{\rm HM}\strut}%
\defop\ECH{{\rm ECH}}%
\def\HFHat{{\widehat{\rm HF}}\strut}%
\defop\Harm{{\cal H}}%
\def\Sobolev#1,#2{\mathop{\rm L}\nolimits_{#1}^{#2}}%
\defop\LTwo{{\rm L}^2}%
\defop\LTwoLoc{{\rm L}^2_{\rm loc}}%
\def\Tan{{\rm T}}%
\defop\Lie{{\rm Lie}}%
\defop\Cliff{{\rm C\ell}}%
\defop\Aut{{\rm Aut}}%
\def\Exterior{{\Lambda}}%
\defop\Ext{{\rm Ext}}%
\defop\Fred{{\rm Fred}}%
\defop\Stab{{\rm Stab}}%
\defop\Pin{{\rm Pin}}%
\def\SpinC{{{\rm Spin}^\C}}%
\def\Spinor{{\cal S}}%
\defop\Interior{{\rm int}}%
\defop\Obj{{\rm ob}}%
\defop\Ball{{\rm B}}%
\defop\supp{{\rm supp}}%
\defop\AssocGrad{{\rm Gr}}%
%
\defop\incl{\hookrightarrow}%
\def\epsilon{\varepsilon}%

\def\Diagram#1{{%
	\def\Up{\Big\uparrow}%
	\def\Down{\Big\downarrow}%
	\def\Right{\longrightarrow}%
	\def\Left{\longleftarrow}%
	\advance\baselineskip by 1em%
	\vcenter{
	\halign{%
		\hfil$##$\hfil&&\quad\hfil$##$\hfil\cr%
		\mathstrut\cr%
		\noalign{\kern-\baselineskip}%
		#1\crcr%
		\mathstrut\cr%
		\noalign{\kern-\baselineskip}%
	}}%
}}%

\defop\SW{{\cal X}}%
\defop\FDSW{{\cal V}}%
\defop\SWF{{\rm SWF}}%
\defop\CSD{{\rm CSD}}%
\def\SpinCStruct{{\frak s}}%
\defop\Inv{{\rm Inv}}%
\defop\Conley{{I}}%
\defop\KohnDirac{{\frak KD}}%
\def\Connex{{\cal A}}%
\def\Disk{{\rm D}}%
\def\ContactThom{{\cal T}}%
\def\B{{\rm B}}%
\def\E{{\rm E}}%
\def\StabHoSpectra#1#2{{\bar{h}#1{\cal S}#2}}%
\defop\SignOp{{\rm Sign}}%

\def\ThomClass{{\theta}}%
\def\BorelHMTilde{\GeomCompl\HMTilde}%
\defop\TB{{\rm tb}}%
\def\Coulomb{{\rm Coul}}%
\def\FreeNumericalInvariant{{F}}%
\def\TorNumericalInvariant{{T}}%
\def\Field{k}%

\global\def\AlexanderBriggsOnTypesCited{0}%
\def\AlexanderBriggsOnTypesAlt{%
	\global\def\AlexanderBriggsOnTypesCited{1}%
	Alexander\ \&\ Briggs\ 1926%
}%
\global\def\AtiyahPatodiSingerSpectralCited{0}%
\global\def\BaragliaHekmatiEquivariantCited{0}%
\def\BaragliaHekmatiEquivariant{%
	\global\def\BaragliaHekmatiEquivariantCited{1}%
	Baraglia\ \&\ Hekmati\ (2021)%
}%
\def\BaragliaHekmatiEquivariantAlt{%
	\global\def\BaragliaHekmatiEquivariantCited{1}%
	Baraglia\ \&\ Hekmati\ 2021%
}%
\global\def\BorelNouvelleCited{0}%
\def\BorelNouvelle{%
	\global\def\BorelNouvelleCited{1}%
	Borel\ (1955)%
}%
\global\def\BorelBredonFloydMontgomeryPalaisSeminarCited{0}%
\global\def\CauchyCoursCited{0}%
\global\def\ChelnokovMednykhHantzscheWendtCited{0}%
\global\def\ColinGhigginiHondaEquivalenceICited{0}%
\global\def\ColinGhigginiHondaEquivalenceIICited{0}%
\global\def\ConleyIsolatedCited{0}%
\def\ConleyIsolated{%
	\global\def\ConleyIsolatedCited{1}%
	Conley\ (1978)%
}%
\def\ConleyIsolatedAlt{%
	\global\def\ConleyIsolatedCited{1}%
	Conley\ 1978%
}%
\global\def\CorneaHomotopicalCited{0}%
\global\def\CostenobleWanerEquivariantPoincareCited{0}%
\global\def\CostenobleWanerEquivariantOrdinaryCited{0}%
\global\def\CristofaroGardinerAbsoluteCited{0}%
\global\def\DonnellyEtaCited{0}%
\global\def\EcheverriaNaturalityCited{0}%
\global\def\EliashbergClassificationCited{0}%
\def\EliashbergClassification{%
	\global\def\EliashbergClassificationCited{1}%
	Eliashberg\ (1989)%
}%
\global\def\FloerRefinementCited{0}%
\global\def\FroyshovMonopoleCited{0}%
\def\FroyshovMonopole{%
	\global\def\FroyshovMonopoleCited{1}%
	Fr\o yshov\ (2010)%
}%
\global\def\FroyshovSeibergWittenCited{0}%
\def\FroyshovSeibergWitten{%
	\global\def\FroyshovSeibergWittenCited{1}%
	Froyshov\ (1996)%
}%
\global\def\GhigginiOzsvathSzaboCited{0}%
\global\def\GhigginiTightCited{0}%
\def\GhigginiTight{%
	\global\def\GhigginiTightCited{1}%
	Ghiggini\ (2008)%
}%
\global\def\GhigginiLiscaStipsiczClassificationCited{0}%
\def\GhigginiLiscaStipsiczClassification{%
	\global\def\GhigginiLiscaStipsiczClassificationCited{1}%
	Ghiggini,\ Lisca\ \&\ Stipsicz\ (2006)%
}%
\global\def\GhigginiHornMorrisTightCited{0}%
\global\def\GirouxStructuresCited{0}%
\def\GirouxStructures{%
	\global\def\GirouxStructuresCited{1}%
	Giroux\ (2000)%
}%
\global\def\GompfHandlebodyCited{0}%
\def\GompfHandlebody{%
	\global\def\GompfHandlebodyCited{1}%
	Gompf\ (1998)%
}%
\global\def\GompfStipsiczFourManifoldsCited{0}%
\global\def\GonzaloBranchedCited{0}%
\def\GonzaloBranched{%
	\global\def\GonzaloBranchedCited{1}%
	Gonzalo\ (1987)%
}%
\global\def\GreenleesMayGeneralizedCited{0}%
\global\def\HendricksLipshitzSarkarFlexibleCited{0}%
\def\HendricksLipshitzSarkarFlexible{%
	\global\def\HendricksLipshitzSarkarFlexibleCited{1}%
	Hendricks,\ Lipshitz\ \&\ Sarkar\ (2016)%
}%
\global\def\HondaClassificationCited{0}%
\def\HondaClassification{%
	\global\def\HondaClassificationCited{1}%
	Honda\ (1999)%
}%
\global\def\HondaClassificationIICited{0}%
\def\HondaClassificationII{%
	\global\def\HondaClassificationIICited{1}%
	Honda\ (2000)%
}%
\global\def\HutchingsTaubesGluingCited{0}%
\global\def\IidaTaniguchiSeibergWittenCited{0}%
\global\def\KangZModTwoEquivariantCited{0}%
\def\KangZModTwoEquivariant{%
	\global\def\KangZModTwoEquivariantCited{1}%
	Kang\ (2018)%
}%
\global\def\KarakurtContactCited{0}%
\def\KarakurtContact{%
	\global\def\KarakurtContactCited{1}%
	Karakurt\ (2014)%
}%
\global\def\KhandhawitNewGaugeCited{0}%
\global\def\KhandhawitLinSasahiraUnfoldedIICited{0}%
\global\def\KhuzamLiftingCited{0}%
\global\def\KronheimerMrowkaMonopolesAndContactStructuresCited{0}%
\def\KronheimerMrowkaMonopolesAndContactStructures{%
	\global\def\KronheimerMrowkaMonopolesAndContactStructuresCited{1}%
	Kronheimer\ \&\ Mrowka\ (1997)%
}%
\def\KronheimerMrowkaMonopolesAndContactStructuresAlt{%
	\global\def\KronheimerMrowkaMonopolesAndContactStructuresCited{1}%
	Kronheimer\ \&\ Mrowka\ 1997%
}%
\global\def\KronheimerMrowkaMonopolesAndThreeManifoldsCited{0}%
\def\KronheimerMrowkaMonopolesAndThreeManifolds{%
	\global\def\KronheimerMrowkaMonopolesAndThreeManifoldsCited{1}%
	Kronheimer\ \&\ Mrowka\ (2007)%
}%
\def\KronheimerMrowkaMonopolesAndThreeManifoldsAlt{%
	\global\def\KronheimerMrowkaMonopolesAndThreeManifoldsCited{1}%
	Kronheimer\ \&\ Mrowka\ 2007%
}%
\global\def\KronheimerMrowkaOzsvathSzaboMonopolesCited{0}%
\def\KronheimerMrowkaOzsvathSzaboMonopoles{%
	\global\def\KronheimerMrowkaOzsvathSzaboMonopolesCited{1}%
	Kronheimer,\ Mrowka,\ Ozsv\'ath\ \&\ Szab\'o\ (2007)%
}%
\global\def\KurlandHomotopyCited{0}%
\global\def\KutluhanLeeTaubesHFEqHMICited{0}%
\global\def\KutluhanLeeTaubesHFEqHMIICited{0}%
\global\def\KutluhanLeeTaubesHFEqHMIIICited{0}%
\global\def\KutluhanLeeTaubesHFEqHMIVCited{0}%
\global\def\KutluhanLeeTaubesHFEqHMVCited{0}%
\global\def\LewisMaySteinbergerEquivariantCited{0}%
\global\def\LidmanManolescuEquivalenceCited{0}%
\def\LidmanManolescuEquivalence{%
	\global\def\LidmanManolescuEquivalenceCited{1}%
	Lidman\ \&\ Manolescu\ (2016)%
}%
\def\LidmanManolescuEquivalenceAlt{%
	\global\def\LidmanManolescuEquivalenceCited{1}%
	Lidman\ \&\ Manolescu\ 2016%
}%
\global\def\LidmanManolescuCoveringCited{0}%
\def\LidmanManolescuCovering{%
	\global\def\LidmanManolescuCoveringCited{1}%
	Lidman\ \&\ Manolescu\ (2018)%
}%
\def\LidmanManolescuCoveringAlt{%
	\global\def\LidmanManolescuCoveringCited{1}%
	Lidman\ \&\ Manolescu\ 2018%
}%
\global\def\LinSolvCited{0}%
\global\def\LinLipnowskiHyperbolicCited{0}%
\def\LinLipnowskiHyperbolic{%
	\global\def\LinLipnowskiHyperbolicCited{1}%
	Lin\ \&\ Lipnowski\ (2022)%
}%
\global\def\LinRubermanSavelievFroyshovCited{0}%
\global\def\LiscaStipsiczOzsvathSzaboICited{0}%
\def\LiscaStipsiczOzsvathSzaboIAlt{%
	\global\def\LiscaStipsiczOzsvathSzaboICited{1}%
	Lisca\ \&\ Stipsicz\ 2004%
}%
\global\def\LiscaStipsiczOzsvathSzaboIICited{0}%
\def\LiscaStipsiczOzsvathSzaboIIYear{%
	\global\def\LiscaStipsiczOzsvathSzaboIICited{1}%
	2007a%
}%
\global\def\LiscaStipsiczOzsvathSzaboIIICited{0}%
\def\LiscaStipsiczOzsvathSzaboIIIYear{%
	\global\def\LiscaStipsiczOzsvathSzaboIIICited{1}%
	2007b%
}%
\global\def\LivingstonSeifertCited{0}%
\def\LivingstonSeifert{%
	\global\def\LivingstonSeifertCited{1}%
	Livingston\ (2002)%
}%
\global\def\ManolescuSeibergWittenFloerCited{0}%
\def\ManolescuSeibergWittenFloer{%
	\global\def\ManolescuSeibergWittenFloerCited{1}%
	Manolescu\ (2003)%
}%
\def\ManolescuSeibergWittenFloerAlt{%
	\global\def\ManolescuSeibergWittenFloerCited{1}%
	Manolescu\ 2003%
}%
\global\def\ManolescuGluingCited{0}%
\global\def\MatkovicClassificationCited{0}%
\def\MatkovicClassification{%
	\global\def\MatkovicClassificationCited{1}%
	Matkovi\v c\ (2018)%
}%
\def\MatkovicClassificationAlt{%
	\global\def\MatkovicClassificationCited{1}%
	Matkovi\v c\ 2018%
}%
\global\def\MayEtAlEquivariantCited{0}%
\def\MayEtAlEquivariantAlt{%
	\global\def\MayEtAlEquivariantCited{1}%
	May\ \&\ al.\ 1996%
}%
\global\def\MischaikowConleyCited{0}%
\global\def\MoserElementaryCited{0}%
\global\def\NicolaescuNotesCited{0}%
\global\def\OzsvathSzaboAbsolutelyCited{0}%
\global\def\OzsvathSzaboHolomorphicCited{0}%
\global\def\OzsvathSzaboContactCited{0}%
\def\OzsvathSzaboContact{%
	\global\def\OzsvathSzaboContactCited{1}%
	Ozsv\'ath\ \&\ Szab\'o\ (2005)%
}%
\global\def\PetitSpinCCited{0}%
\def\PetitSpinCAlt{%
	\global\def\PetitSpinCCited{1}%
	Petit\ 2005%
}%
\global\def\PlamenevskayaContactStructuresWithDistinctCited{0}%
\def\PlamenevskayaContactStructuresWithDistinctAlt{%
	\global\def\PlamenevskayaContactStructuresWithDistinctCited{1}%
	Plamenevskaya\ 2004%
}%
\global\def\RamosAbsoluteCited{0}%
\global\def\RolfsenKnotsCited{0}%
\def\RolfsenKnots{%
	\global\def\RolfsenKnotsCited{1}%
	Rolfsen\ (2003)%
}%
\global\def\RosoSeibergWittenCited{0}%
\def\RosoSeibergWitten{%
	\global\def\RosoSeibergWittenCited{1}%
	Roso\ (2023)%
}%
\def\RosoSeibergWittenAlt{%
	\global\def\RosoSeibergWittenCited{1}%
	Roso\ 2023%
}%
\global\def\ShubinPseudodifferentialCited{0}%
\global\def\SmithTransformationsCited{0}%
\def\SmithTransformations{%
	\global\def\SmithTransformationsCited{1}%
	Smith\ (1938)%
}%
\global\def\TaubesSeibergWittenSymplecticCited{0}%
\global\def\TaubesSWToGrCited{0}%
\global\def\TaubesWeinsteinICited{0}%
\def\TaubesWeinsteinIAlt{%
	\global\def\TaubesWeinsteinICited{1}%
	Taubes\ 2007%
}%
\global\def\TaubesWeinsteinIICited{0}%
\def\TaubesWeinsteinII{%
	\global\def\TaubesWeinsteinIICited{1}%
	Taubes\ (2009)%
}%
\global\def\TaubesEmbeddedICited{0}%
\global\def\TaubesEmbeddedIICited{0}%
\global\def\TaubesEmbeddedIIICited{0}%
\global\def\TaubesEmbeddedIVCited{0}%
\global\def\TaubesEmbeddedVCited{0}%
\global\def\TomDieckTransformationCited{0}%
\global\def\WuLegendrianCited{0}%
\def\WuLegendrian{%
	\global\def\WuLegendrianCited{1}%
	Wu\ (2006)%
}%
\def\Bib{{%
	\parindent=0cm%
	\if\AlexanderBriggsOnTypesCited 1%
		{\bf Alexander,\ J.W.\ \&\ Briggs,\ G.B.}
		(1926).
		``On types of knotted curves.''
		{\it Annals of Mathematics},
		pp. 562--586.
		JSTOR.
		\par%
	\fi%
	\if\AtiyahPatodiSingerSpectralCited 1%
		{\bf Atiyah,\ M.F.,\ Patodi,\ V.K.\ \&\ Singer,\ I.M.}
		(1975).
		``Spectral Asymmetry and Riemannian Geometry I.''
		{\it Mathematical Proceedings of the Cambridge Philosophical Society},
		77(1),
		pp. 43--69.
		Cambridge University Press.
		\par%
	\fi%
	\if\BaragliaHekmatiEquivariantCited 1%
		{\bf Baraglia,\ D.\ \&\ Hekmati,\ P.}
		(2021).
		``Equivariant Seiberg--Witten--Floer cohomology.''
		{\it arXiv preprint arXiv:2108.06855}.
		\par%
	\fi%
	\if\BorelNouvelleCited 1%
		{\bf Borel,\ A.}
		(1955).
		``Nouvelle d{\'e}monstration d'un th{\'e}oreme de PA Smith.''
		{\it Comment. Math. Helv},
		29,
		pp. 27--39.
		\par%
	\fi%
	\if\BorelBredonFloydMontgomeryPalaisSeminarCited 1%
		{\bf Borel,\ A.,\ Bredon,\ G.,\ Floyd,\ E.E.,\ Montgomery,\ D.\ \&\ Palais,\ R.}
		(1960).
		{\it Seminar on Transformation Groups. (AM-46)}.
		Princeton University Press.
		\par%
	\fi%
	\if\CauchyCoursCited 1%
		{\bf Cauchy,\ A.L.}
		(1821).
		{\it Cours d'Analyse de l'\'Ecole Royale Polytechnique}.
		de Bure, Paris.
		\par%
	\fi%
	\if\ChelnokovMednykhHantzscheWendtCited 1%
		{\bf Chelnokov,\ G.\ \&\ Mednykh,\ A.}
		(2020).
		{\it On the coverings of Hantzsche-Wendt manifold}.
		arXiv preprint arXiv:2009.06691.
		\par%
	\fi%
	\if\ColinGhigginiHondaEquivalenceICited 1%
		{\bf Colin,\ V.,\ Ghiggini,\ P.\ \&\ Honda,\ K.}
		(2012a).
		{\it The equivalence of Heegaard Floer and embedded contact homology via open book decompositions I}.
		arXiv preprint arXiv:1208.1074.
		\par%
	\fi%
	\if\ColinGhigginiHondaEquivalenceIICited 1%
		{\bf Colin,\ V.,\ Ghiggini,\ P.\ \&\ Honda,\ K.}
		(2012b).
		{\it The equivalence of Heegaard Floer and embedded contact homology via open book decompositions II}.
		arXiv preprint arXiv:1208.1077.
		\par%
	\fi%
	\if\ConleyIsolatedCited 1%
		{\bf Conley,\ C.}
		(1978).
		{\it Isolated Invariant Sets and the Morse Index},
		(38).
		American Mathematical Society.
		\par%
	\fi%
	\if\CorneaHomotopicalCited 1%
		{\bf Cornea,\ O.}
		(2000).
		``Homotopical dynamics: suspension and duality.''
		{\it Ergodic Theory and Dynamical Systems},
		20(2),
		pp. 379--391.
		Cambridge University Press.
		\par%
	\fi%
	\if\CostenobleWanerEquivariantPoincareCited 1%
		{\bf Costenoble,\ S.R.\ \&\ Waner,\ S.}
		(1992).
		``Equivariant Poincar\'e Duality.''
		{\it Michigan Mathematical Journal},
		39(2),
		pp. 325--351.
		\par%
	\fi%
	\if\CostenobleWanerEquivariantOrdinaryCited 1%
		{\bf Costenoble,\ S.R.\ \&\ Waner,\ S.}
		(2016).
		{\it Equivariant ordinary homology and cohomology}.
		Springer.
		\par%
	\fi%
	\if\CristofaroGardinerAbsoluteCited 1%
		{\bf Cristofaro-Gardiner,\ D.}
		(2013).
		``The absolute gradings on embedded contact homology and Seiberg--Witten Floer cohomology.''
		{\it Algebraic \& Geometric Topology},
		13(4),
		pp. 2239--2260.
		Mathematical Sciences Publishers.
		\par%
	\fi%
	\if\DonnellyEtaCited 1%
		{\bf Donnelly,\ H.}
		(1978).
		``Eta invariants for $G$-spaces.''
		{\it Indiana University Mathematics Journal},
		27(6),
		pp. 889--918.
		\par%
	\fi%
	\if\EcheverriaNaturalityCited 1%
		{\bf Echeverria,\ M.}
		(2020).
		``Naturality of the contact invariant in monopole Floer homology under strong symplectic cobordisms.''
		{\it Algebraic \& Geometric Topology},
		20(4).
		\par%
	\fi%
	\if\EliashbergClassificationCited 1%
		{\bf Eliashberg,\ Ya.}
		(1989).
		``Classification of overtwisted contact structures on 3-manifolds.''
		{\it Inventiones mathematicae},
		98(3),
		pp. 623--637.
		Citeseer.
		\par%
	\fi%
	\if\FloerRefinementCited 1%
		{\bf Floer,\ A.}
		(1987).
		``A Refinement of the Conley Index and an Application to the Stability of Hyperbolic Invariant Sets.''
		{\it Ergodic Theory and Dynamical Systems},
		7(1).
		\par%
	\fi%
	\if\FroyshovMonopoleCited 1%
		{\bf Fr\o yshov,\ K.A.}
		(2010).
		``Monopole Floer homology for rational homology 3-spheres.''
		{\it Duke Mathematical Journal},
		155(3),
		pp. 519--576.
		Duke University Press.
		\par%
	\fi%
	\if\FroyshovSeibergWittenCited 1%
		{\bf Froyshov,\ K.A.}
		(1996).
		``The Seiberg--Witten equations and four-manifolds with boundary.''
		{\it Mathematical Research Letters},
		3(3),
		pp. 373--390.
		International Press of Boston.
		\par%
	\fi%
	\if\GhigginiOzsvathSzaboCited 1%
		{\bf Ghiggini,\ P.}
		(2006).
		``Ozsv\'ath-Szab\'o invariants and fillability of contact structures.''
		{\it Mathematische Zeitschrift},
		253(1).
		\par%
	\fi%
	\if\GhigginiTightCited 1%
		{\bf Ghiggini,\ P.}
		(2008).
		``On tight contact structures with negative maximal twisting number on small Seifert manifolds.''
		{\it Algebraic \& Geometric Topology},
		8(1),
		pp. 381--396.
		\par%
	\fi%
	\if\GhigginiLiscaStipsiczClassificationCited 1%
		{\bf Ghiggini,\ P.,\ Lisca,\ P.\ \&\ Stipsicz,\ A.I.}
		(2006).
		``Classification of Tight Contact Structures on Small Seifert 3-Manifolds with $e_0\geq0$.''
		{\it Proceedings of the American Mathematical Society},
		pp. 909--916.
		\par%
	\fi%
	\if\GhigginiHornMorrisTightCited 1%
		{\bf Ghiggini,\ P.\ \&\ Van Horn-Morris,\ J.}
		(2016).
		``Tight contact structures on the Brieskorn spheres $-\Sigma(2,3,6n-1)$ and contact invariants.''
		{\it Journal f\"ur die reine und angewandte Mathematik},
		2016(718),
		pp. 1--24.
		\par%
	\fi%
	\if\GirouxStructuresCited 1%
		{\bf Giroux,\ E.}
		(2000).
		``Structures de contact en dimension trois et bifurcations des feuilletages de surfaces.''
		{\it Inventiones Mathematicae},
		141(3),
		pp. 615--689.
		\par%
	\fi%
	\if\GompfHandlebodyCited 1%
		{\bf Gompf,\ R.E.}
		(1998).
		``Handlebody construction of Stein surfaces.''
		{\it Annals of mathematics},
		pp. 619--693.
		\par%
	\fi%
	\if\GompfStipsiczFourManifoldsCited 1%
		{\bf Gompf,\ R.E.\ \&\ Stipsicz,\ A.I.}
		(1999).
		{\it 4-manifolds and Kirby calculus},
		(20).
		American Mathematical Soc..
		\par%
	\fi%
	\if\GonzaloBranchedCited 1%
		{\bf Gonzalo,\ J.}
		(1987).
		``Branched covers and contact structures.''
		{\it Proceedings of the American Mathematical Society},
		101(2),
		pp. 347--352.
		\par%
	\fi%
	\if\GreenleesMayGeneralizedCited 1%
		{\bf Greenlees,\ J.P.C.\ \&\ May,\ J.P.}
		(1995).
		{\it Generalized Tate cohomology},
		543.
		American Mathematical Soc..
		\par%
	\fi%
	\if\HendricksLipshitzSarkarFlexibleCited 1%
		{\bf Hendricks,\ K.,\ Lipshitz,\ R.\ \&\ Sarkar,\ S.}
		(2016).
		``A flexible construction of equivariant Floer homology and applications.''
		{\it Journal of Topology},
		9(4),
		pp. 1153--1236.
		\par%
	\fi%
	\if\HondaClassificationCited 1%
		{\bf Honda,\ K.}
		(1999).
		``On the classification of tight contact structures I.''
		{\it Geometry \& Topology},
		4,
		pp. 309--368.
		\par%
	\fi%
	\if\HondaClassificationIICited 1%
		{\bf Honda,\ K.}
		(2000).
		``On the classification of tight contact structures II.''
		{\it Journal of Differential Geometry},
		55(1),
		pp. 83--143.
		Lehigh University.
		\par%
	\fi%
	\if\HutchingsTaubesGluingCited 1%
		{\bf Hutchings,\ M.\ \&\ Taubes,\ C.H.}
		(2007).
		``Gluing pseudoholomorphic curves along branched covered cylinders I.''
		{\it Journal of Symplectic Geometry},
		5(1),
		pp. 43--137.
		International Press of Boston.
		\par%
	\fi%
	\if\IidaTaniguchiSeibergWittenCited 1%
		{\bf Iida,\ N.\ \&\ Taniguchi,\ M.}
		(2021).
		``Seiberg--Witten Floer Homotopy Contact Invariant.''
		{\it Studia Scientiarum Mathematicarum Hungarica},
		58(4),
		pp. 505--558.
		Akad\'emiai Kiad\'o Budapest.
		\par%
	\fi%
	\if\KangZModTwoEquivariantCited 1%
		{\bf Kang,\ S.}
		(2018).
		{\it $Z_2$-equivariant Heegaard Floer cohomology and transverse knots}.
		Ph.D.\ Dissertation, University of Oxford.
		\par%
	\fi%
	\if\KarakurtContactCited 1%
		{\bf Karakurt,\ \c C.}
		(2014).
		``Contact structures on plumbed 3-manifolds.''
		{\it Kyoto Journal of Mathematics},
		54(2),
		pp. 271--294.
		Duke University Press.
		\par%
	\fi%
	\if\KhandhawitNewGaugeCited 1%
		{\bf Khandhawit,\ T.}
		(2015).
		``A new gauge slice for the relative Bauer--Furuta invariants.''
		{\it Geometry \& Topology},
		19(3),
		pp. 1631--1655.
		Mathematical Sciences Publishers.
		\par%
	\fi%
	\if\KhandhawitLinSasahiraUnfoldedIICited 1%
		{\bf Khandhawit,\ T.,\ Lin,\ J.\ \&\ Sasahira,\ H.}
		(2018).
		{\it Unfolded Seiberg--Witten Floer spectra, II: Relative invariants and the gluing theorem}.
		arXiv preprint arXiv:1809.09151.
		\par%
	\fi%
	\if\KhuzamLiftingCited 1%
		{\bf Khuzam,\ M.B.}
		(2012).
		``Lifting the $3$-dimensional invariant of $2$-plane fields on $3$-manifolds.''
		{\it Topology and its Applications},
		159(3),
		pp. 704--710.
		\par%
	\fi%
	\if\KronheimerMrowkaMonopolesAndContactStructuresCited 1%
		{\bf Kronheimer,\ P.\ \&\ Mrowka,\ T.}
		(1997).
		``Monopoles and Contact Structures.''
		{\it Inventiones Mathematicae},
		130(2),
		pp. 209--255.
		\par%
	\fi%
	\if\KronheimerMrowkaMonopolesAndThreeManifoldsCited 1%
		{\bf Kronheimer,\ P.\ \&\ Mrowka,\ T.}
		(2007).
		{\it Monopoles and three-manifolds},
		10.
		Cambridge University Press.
		\par%
	\fi%
	\if\KronheimerMrowkaOzsvathSzaboMonopolesCited 1%
		{\bf Kronheimer,\ P.,\ Mrowka,\ T.,\ Ozsv\'ath,\ P.\ \&\ Szab\'o,\ Z.}
		(2007).
		``Monopoles and Lens Space Surgeries.''
		{\it Annals of Mathematics},
		pp. 457--546.
		JSTOR.
		\par%
	\fi%
	\if\KurlandHomotopyCited 1%
		{\bf Kurland,\ H.L.}
		(1982).
		``Homotopy Invariants of Repeller-Attractor Pairs. I. The P\"uppe Sequence of an RA Pair.''
		{\it Journal of Differential Equations},
		46(1),
		pp. 1--31.
		\par%
	\fi%
	\if\KutluhanLeeTaubesHFEqHMICited 1%
		{\bf Kutluhan,\ \c C.,\ Lee,\ Y.J.\ \&\ Taubes,\ C.H.}
		(2020a).
		``HF=HM, I: Heegaard Floer homology and Seiberg--Witten Floer homology.''
		{\it Geometry \& Topology},
		24(6),
		pp. 2829--2854.
		Mathematical Sciences Publishers.
		\par%
	\fi%
	\if\KutluhanLeeTaubesHFEqHMIICited 1%
		{\bf Kutluhan,\ \c C.,\ Lee,\ Y.J.\ \&\ Taubes,\ C.H.}
		(2020b).
		``HF=HM, II: Reeb orbits and holomorphic curves for the ech/Heegaard Floer correspondence.''
		{\it Geometry \& Topology},
		24(6),
		pp. 2855--3012.
		Mathematical Sciences Publishers.
		\par%
	\fi%
	\if\KutluhanLeeTaubesHFEqHMIIICited 1%
		{\bf Kutluhan,\ \c C.,\ Lee,\ Y.J.\ \&\ Taubes,\ C.H.}
		(2020c).
		``HF=HM, III: holomorphic curves and the differential for the ech/Heegaard Floer correspondenc.''
		{\it Geometry \& Topology},
		24(6),
		pp. 3013--3218.
		Mathematical Sciences Publishers.
		\par%
	\fi%
	\if\KutluhanLeeTaubesHFEqHMIVCited 1%
		{\bf Kutluhan,\ \c C.,\ Lee,\ Y.J.\ \&\ Taubes,\ C.H.}
		(2021a).
		``HF=HM, IV: The Seiberg--Witten Floer homology and ech correspondence.''
		{\it Geometry \& Topology},
		24(7),
		pp. 3219--3469.
		Mathematical Sciences Publishers.
		\par%
	\fi%
	\if\KutluhanLeeTaubesHFEqHMVCited 1%
		{\bf Kutluhan,\ \c C.,\ Lee,\ Y.J.\ \&\ Taubes,\ C.H.}
		(2021b).
		``HF=HM, V: Seiberg--Witten Floer homology and handle additions.''
		{\it Geometry \& Topology},
		24(7),
		pp. 3471--3748.
		Mathematical Sciences Publishers.
		\par%
	\fi%
	\if\LewisMaySteinbergerEquivariantCited 1%
		{\bf Lewis Jr.,\ L.G.,\ May,\ J.P.\ \&\ Steinberger,\ M.}
		(1986).
		``Equivariant stable homotopy theory.''
		{\it Lecture notes in mathematics},
		1213.
		Springer-Verlag.
		\par%
	\fi%
	\if\LidmanManolescuEquivalenceCited 1%
		{\bf Lidman,\ T.\ \&\ Manolescu,\ C.}
		(2016).
		{\it The equivalence of two Seiberg--Witten Floer homologies}.
		arXiv preprint arXiv:1603.00582.
		\par%
	\fi%
	\if\LidmanManolescuCoveringCited 1%
		{\bf Lidman,\ T.\ \&\ Manolescu,\ C.}
		(2018).
		``Floer homology and covering spaces.''
		{\it Geometry \& Topology},
		22(5),
		pp. 2817--2838.
		\par%
	\fi%
	\if\LinSolvCited 1%
		{\bf Lin,\ F.}
		(2020).
		``Monopole Floer homology and SOLV geometry.''
		{\it Annales Henri Lebesgue},
		3,
		pp. 1117--1131.
		\par%
	\fi%
	\if\LinLipnowskiHyperbolicCited 1%
		{\bf Lin,\ F.\ \&\ Lipnowski,\ M.}
		(2022).
		``The Seiberg--Witten equations and the length spectrum of hyperbolic three-manifolds.''
		{\it Journal of the American Mathematical Society},
		35(1),
		pp. 233--293.
		\par%
	\fi%
	\if\LinRubermanSavelievFroyshovCited 1%
		{\bf Lin,\ J.,\ Ruberman,\ D.\ \&\ Saveliev,\ N.}
		(2018).
		{\it On the Fr\o yshov invariant and monopole Lefschetz number}.
		arXiv preprint arXiv:1802.07704.
		\par%
	\fi%
	\if\LiscaStipsiczOzsvathSzaboICited 1%
		{\bf Lisca,\ P.\ \&\ Stipsicz,\ A.I.}
		(2004).
		``Ozsv{\'a}th--Szab{\'o} invariants and tight contact three-manifolds I.''
		{\it Geometry \& Topology},
		8(2),
		pp. 925--945.
		Mathematical Sciences Publishers.
		\par%
	\fi%
	\if\LiscaStipsiczOzsvathSzaboIICited 1%
		{\bf Lisca,\ P.\ \&\ Stipsicz,\ A.I.}
		(2007a).
		``Ozsv{\'a}th-Szab{\'o} invariants and tight contact three-manifolds, II.''
		{\it Journal of Differential Geometry},
		75(1).
		\par%
	\fi%
	\if\LiscaStipsiczOzsvathSzaboIIICited 1%
		{\bf Lisca,\ P.\ \&\ Stipsicz,\ A.I.}
		(2007b).
		``Ozsv{\'a}th-Szab{\'o} invariants and tight contact 3-manifolds, III.''
		{\it Journal of Symplectic Geometry},
		5(4),
		pp. 357--384.
		International Press of Boston.
		\par%
	\fi%
	\if\LivingstonSeifertCited 1%
		{\bf Livingston,\ C.}
		(2002).
		``Seifert forms and concordance.''
		{\it Geometry \& Topology},
		6(1),
		pp. 403--408.
		Mathematical Sciences Publishers.
		\par%
	\fi%
	\if\ManolescuSeibergWittenFloerCited 1%
		{\bf Manolescu,\ C.}
		(2003).
		``Seiberg--Witten--Floer stable homotopy type of three-manifolds with $b_1=0$.''
		{\it Geometry \& Topology},
		7(2),
		pp. 889--932.
		Mathematical Sciences Publishers.
		\par%
	\fi%
	\if\ManolescuGluingCited 1%
		{\bf Manolescu,\ C.}
		(2007).
		``A gluing theorem for the relative Bauer-Furuta invariants.''
		{\it Journal of Differential Geometry},
		76(1),
		pp. 117--153.
		\par%
	\fi%
	\if\MatkovicClassificationCited 1%
		{\bf Matkovi\v c,\ I.}
		(2018).
		``Classification of tight contact structures on small Seifert fibered L–spaces.''
		{\it Algebraic \& Geometric Topology},
		18(1),
		pp. 111--152.
		\par%
	\fi%
	\if\MayEtAlEquivariantCited 1%
		{\bf May,\ J.P.\ \&\ al.}
		(1996).
		{\it Equivariant Homotopy and Cohomology Theory},
		(91).
		American Mathematical Soc..
		\par%
	\fi%
	\if\MischaikowConleyCited 1%
		{\bf Mischaikow,\ K.}
		(1995).
		``Conley index theory.''
		{\it Dynamical Systems},
		pp. 119--207.
		Springer.
		\par%
	\fi%
	\if\MoserElementaryCited 1%
		{\bf Moser,\ L.}
		(1971).
		``Elementary surgery along a torus knot.''
		{\it Pacific Journal of Mathematics},
		38(3),
		pp. 737--745.
		\par%
	\fi%
	\if\NicolaescuNotesCited 1%
		{\bf Nicolaescu,\ L.I.}
		(2000).
		{\it Notes on Seiberg--Witten Theory},
		28.
		American Mathematical Soc..
		\par%
	\fi%
	\if\OzsvathSzaboAbsolutelyCited 1%
		{\bf Ozsv\'ath,\ P.\ \&\ Szab\'o,\ Z.}
		(2003).
		``Absolutely graded Floer homologies and intersection forms for four-manifolds with boundary.''
		{\it Advances in Mathematics},
		173(2),
		pp. 179--261.
		Elsevier.
		\par%
	\fi%
	\if\OzsvathSzaboHolomorphicCited 1%
		{\bf Ozsv\'ath,\ P.\ \&\ Szab\'o,\ Z.}
		(2004).
		``Holomorphic disks and topological invariants for closed three-manifolds.''
		{\it Annals of Mathematics},
		pp. 1027--1158.
		JSTOR.
		\par%
	\fi%
	\if\OzsvathSzaboContactCited 1%
		{\bf Ozsv\'ath,\ P.\ \&\ Szab\'o,\ Z.}
		(2005).
		``Heegaard Floer homology and contact structures.''
		{\it Duke Mathematical Journal},
		129(1).
		\par%
	\fi%
	\if\PetitSpinCCited 1%
		{\bf Petit,\ R.}
		(2005).
		``Spinc-structures and Dirac operators on contact manifolds.''
		{\it Differential Geometry and its Applications},
		22(2),
		pp. 229--252.
		Elsevier.
		\par%
	\fi%
	\if\PlamenevskayaContactStructuresWithDistinctCited 1%
		{\bf Plamenevskaya,\ O.}
		(2004).
		``Contact structures with distinct Heegaard Floer invariants.''
		{\it Mathematical Research Letters},
		11,
		pp. 547--561.
		\par%
	\fi%
	\if\RamosAbsoluteCited 1%
		{\bf Ramos,\ V.G.B.}
		(2018).
		``Absolute gradings on ECH and Heegaard Floer homology.''
		{\it Quantum Topology},
		9(2),
		pp. 207--228.
		\par%
	\fi%
	\if\RolfsenKnotsCited 1%
		{\bf Rolfsen,\ D.}
		(2003).
		{\it Knots and links},
		346.
		American Mathematical Soc..
		\par%
	\fi%
	\if\RosoSeibergWittenCited 1%
		{\bf Roso,\ B.R.S.}
		(2023).
		``Seiberg--Witten Floer spectra and contact structures.''
		{\it Journal of Fixed Point Theory and Applications},
		25(2).
		Springer.
		\par%
	\fi%
	\if\ShubinPseudodifferentialCited 1%
		{\bf Shubin,\ M.A.}
		(2000).
		{\it Pseudodifferential Operators and Spectral Theory},
		2.
		Springer.
		\par%
	\fi%
	\if\SmithTransformationsCited 1%
		{\bf Smith,\ P.A.}
		(1938).
		``Transformations of finite period.''
		{\it Annals of Mathematics},
		pp. 127--164.
		JSTOR.
		\par%
	\fi%
	\if\TaubesSeibergWittenSymplecticCited 1%
		{\bf Taubes,\ C.H.}
		(1994).
		``The Seiberg--Witten invariants and symplectic forms.''
		{\it Mathematical Research Letters},
		1(6),
		pp. 809--822.
		International Press of Boston.
		\par%
	\fi%
	\if\TaubesSWToGrCited 1%
		{\bf Taubes,\ C.H.}
		(1996).
		``SW $\Rightarrow$ Gr: from the Seiberg--Witten equations to pseudo-holomorphic curves.''
		{\it Journal of the American Mathematical Society},
		pp. 845--918.
		JSTOR.
		\par%
	\fi%
	\if\TaubesWeinsteinICited 1%
		{\bf Taubes,\ C.H.}
		(2007).
		``The Seiberg--Witten equations and the Weinstein conjecture.''
		{\it Geometry \& Topology},
		11(4),
		pp. 2117--2202.
		Mathematical Sciences Publishers.
		\par%
	\fi%
	\if\TaubesWeinsteinIICited 1%
		{\bf Taubes,\ C.H.}
		(2009).
		``The Seiberg--Witten equations and the Weinstein conjecture II: More closed integral curves of the Reeb vector field.''
		{\it Geometry \& Topology},
		13(3),
		pp. 1337--1417.
		Mathematical Sciences Publishers.
		\par%
	\fi%
	\if\TaubesEmbeddedICited 1%
		{\bf Taubes,\ C.H.}
		(2010a).
		``Embedded contact homology and Seiberg--Witten Floer cohomology I.''
		{\it Geometry \& Topology},
		14(5),
		pp. 2497--2581.
		Mathematical Sciences Publishers.
		\par%
	\fi%
	\if\TaubesEmbeddedIICited 1%
		{\bf Taubes,\ C.H.}
		(2010b).
		``Embedded contact homology and Seiberg--Witten Floer cohomology II.''
		{\it Geometry \& Topology},
		14(5),
		pp. 2583--2720.
		Mathematical Sciences Publishers.
		\par%
	\fi%
	\if\TaubesEmbeddedIIICited 1%
		{\bf Taubes,\ C.H.}
		(2010c).
		``Embedded contact homology and Seiberg--Witten Floer cohomology III.''
		{\it Geometry \& Topology},
		14(5),
		pp. 2721--2817.
		Mathematical Sciences Publishers.
		\par%
	\fi%
	\if\TaubesEmbeddedIVCited 1%
		{\bf Taubes,\ C.H.}
		(2010d).
		``Embedded contact homology and Seiberg--Witten Floer cohomology IV.''
		{\it Geometry \& Topology},
		14(5),
		pp. 2819--2960.
		Mathematical Sciences Publishers.
		\par%
	\fi%
	\if\TaubesEmbeddedVCited 1%
		{\bf Taubes,\ C.H.}
		(2010e).
		``Embedded contact homology and Seiberg--Witten Floer cohomology V.''
		{\it Geometry \& Topology},
		14(5),
		pp. 2961--3000.
		Mathematical Sciences Publishers.
		\par%
	\fi%
	\if\TomDieckTransformationCited 1%
		{\bf tom Dieck,\ T.}
		(1987).
		{\it Transformation Groups}.
		W. de Gruyter.
		\par%
	\fi%
	\if\WuLegendrianCited 1%
		{\bf Wu,\ H.}
		(2006).
		``Legendrian vertical circles in small Seifert spaces.''
		{\it Communications in Contemporary Mathematics},
		8(2),
		pp. 219--246.
		\par%
	\fi%
}}%

\Title{On contact modulo p L-space covers}
\Abs
The purpose of this {\Document}
is to extend certain results
of {\RosoSeibergWitten}
which concerned equivariant contact structures
on {\it minimal L-spaces\/}
to the more general setting of
{\it mod $p$ L-spaces}.
This is achieved by considering
the Serre spectral sequence
of the Borel Floer cohomology
as showcased in \BaragliaHekmatiEquivariant.
Along the way,
the author introduces two new numerical invariants
of equivariant contact structures
which emerge from the structure
of the Borel Floer cohomology
as a module over the group cohomology.
\Keywords{
	Contact structures,
	Seiberg--Witten Floer theory,
	regular coverings.
}
\MSC{
	57K33,
	57R58,
	57M10.
}

\Sect{Introduction}
\par
Floer theoretic invariants of contact structures
have enjoyed great success
in determining tightness of contact structures
(vid.\ e.g.\ \LiscaStipsiczOzsvathSzaboIAlt,
{\LiscaStipsiczOzsvathSzaboIIYear},
{\LiscaStipsiczOzsvathSzaboIIIYear})
and also in distinguishing between tight contact structures
where the underlying hyperplane fields are homotopic
(vid.\ e.g.\ \PlamenevskayaContactStructuresWithDistinctAlt).
The primary goal of the present {\Document}
is to elucidate the relationship
between the contact invariants 
of a contact manifold
and its regular elementary abelian covers.
Precise answers as to vanishing and non-vanishing
of the invariants involved shall be given in the case of
mod $p$ L-spaces,
whereas,
in more generality,
a pair of new numerical invariants shall be introduced
which shed some light on the lifting behaviour
of the contact invariant.
The technical apparatus which shall be invoked
in this pursuit shall be that of {\it equivariant\/} Floer theory
as studied by {\BaragliaHekmatiEquivariant}
and the equivariant contact invariant first introduced
in \RosoSeibergWitten,
wherein some of the results presented here 
were presaged in much less generality.
\par
Suppose $(Y,\lambda)$
denote a contact rational homology $3$-sphere.
To the contact form $\lambda$,
there is a naturally associated $\SpinC$ structure
$\SpinCStruct_\lambda$.
In the context of Seiberg--Witten Floer theory,
the contact element
$\psi(\lambda)\in\HMTilde^*(Y,\SpinCStruct_\lambda)$
was first precisely defined in the works of
{\KronheimerMrowkaOzsvathSzaboMonopoles}
and {\TaubesWeinsteinII},
though the fundamentals had essentially already appeared
in \KronheimerMrowkaMonopolesAndContactStructures.
\par
Considering the scenario in which there be
a finite group $G$ acting freely on $Y$
and the contact form $\lambda$
be $G$-equivariant,
it is interesting to relate
the invariant $\psi(\lambda)$
to the invariant $\psi(\lambda/G)$
of the quotient contact manifold
$(Y/G,\lambda/G)$.
In this context,
the following shall be proven in the present \Document.
\Thm{
	\Label\ModPSpinCLSpaceLiftContactInvariant
	Let $G$ be an elementary abelian $p$-group
	for $p$ prime
	and let $\Field$ denote the field of $p$ elements.
	Suppose that $G$ act freely 
	on a rational homology $3$-sphere $Y$
	while preserving a contact form $\lambda$.
	Assume further that
	$\dim_{\Field}\HMTilde^*(Y,\SpinCStruct_\lambda;\Field)=1$.
	It follows that,
	if $\psi(\lambda/G)=0$,
	then $\psi(\lambda)=0$.
	Conversely,
	if $\psi(\lambda/G)\neq0$,
	then
	$\psi(\lambda)\neq0$
	if and only if
	$d_3(\lambda)+1/2=-2h(\SpinCStruct_\lambda)$,
	where $d_3$ denotes
	the $3$-dimensional hyperplane field invariant
	of {\GompfHandlebody}
	and $h$
	denotes the invariant of {\FroyshovMonopole}.
}
\Rem{
	The contact invariants
	and the Fr\o yshov invariant
	here
	are defined with coefficients
	in $\Field$.
}
The proof of this result
shall rely
on the $G$-equivariant
contact class
$\psi_G(\lambda)$
living in the
$G$-equivariant Borel mono\-pole Floer cohomology
$\BorelHMTilde_G^*(Y,\SpinCStruct_\lambda)$.
In \RosoSeibergWitten,
Theorems 7.31 and 7.34,
the author proved a similar result to the above
but concerning only the case when $Y$ be a
{\it minimal L-space}.
{\LinLipnowskiHyperbolic}
introduced the nomenclature minimal L-space
to refer to a rational homology $3$-sphere
admitting a metric for which there be
no irreducible solution to the Seiberg--Witten equations.
Hence,
the main purpose of this article
is to relax the stringent condition of minimality
to the purely cohomological requirement seen above.
In order to achieve this,
the key ingredient shall be
the Serre spectral sequence.
\Thm{
	\Label\SerreSpectralSequence
	(\BaragliaHekmatiEquivariantAlt, Theorem 1.1)
	For $(Y,\SpinCStruct)$
	a $\SpinC$ rational homology $3$-sphere
	with an action a finite group $G$,
	there exists a spectral sequence
	converging to the Borel Floer cohomology
	$\BorelHMTilde_G^*(Y,\SpinCStruct)$
	whose second page is
	$$
		E^{i,j}_2
		=\Homology^i(G,\HMTilde^j(Y,\SpinCStruct)),
	$$
	where $\Homology^*(G,M)$
	is the group cohomology of $G$
	with respect to a $\Z[G]$-module $M$.
}
One manner of ensuring that the condition
$\dim_{\Field}\HMTilde^*(Y,\SpinCStruct_\lambda;\Field)=1$
hold
is to require that $Y$ be a mod $p$ L-space.
\Def{
	A rational homology $3$-sphere $Y$
	is called a {\it mod $p$ L-space\/}
	when,
	for any $\SpinC$ structure $\SpinCStruct$,
	$\dim_{\Field}\HMTilde^*(Y,\SpinCStruct;\Field)=1$.
}
A more popular definition is the following.
\Def{
	A rational homology $3$-sphere $Y$
	is called an {\it L-space\/}
	when,
	for any $\SpinC$ structure $\SpinCStruct$,
	$\HMTilde^*(Y,\SpinCStruct;\Z)=\Z$.
}
The only currently known source of mod $p$ L-spaces
are L-spaces themselves,
which,
by universal coefficients,
are mod $p$ L-spaces
for all $p$.
This is because there are no currently known examples of torsion
in the monopole Floer cohomology of rational homology spheres.
\Cor{
	\Label\ModPLSpaceLiftContactInvariant
	Let $G$ be an elementary abelian $p$-group
	for $p$ prime.
	Suppose that $G$ act freely 
	on a mod $p$ L-space $Y$
	while preserving some contact form $\lambda$.
	If $\psi(\lambda/G)=0$,
	it follows that $\psi(\lambda)=0$.
	Conversely,
	if $\psi(\lambda/G)\neq0$,
	then
	$\psi(\lambda)\neq0$
	if and only if
	$d_3(\lambda)+1/2=-2h(\SpinCStruct_\lambda)$.
}
A well known phenomenon in contact topology
is that a tight contact manifold $(Y/G,\lambda/G)$
may be covered by an overtwisted contact manifold
$(Y,\lambda)$.
In this case,
the contact structure $\lambda/G$ is called
{\it virtually overtwisted.}
Given a tight $(Y/G,\lambda/G)$,
determining whether $(Y,\lambda)$ is tight or overtwisted
is typically a difficult problem
solved only in very special scenarios
such as coverings of
thickened tori,
solid tori
and
lens spaces;
vid.\ \HondaClassificationII, \S1.2.
\par
The usefulness of {\ModPLSpaceLiftContactInvariant}
or,
more generally,
{\ModPSpinCLSpaceLiftContactInvariant}
becomes evident
in a situation
where one desires to determine tightness
of the contact structure $\lambda$
when all one knows is that $\lambda/G$ has non-vanishing contact class.
By performing a Kirby calculus procedure
explained in {\RosoSeibergWitten}, \S8,
one can often compute the value of $d_3(\lambda)$
from that of $d_3(\lambda/G)$
and,
thereby,
determine whether $\psi(\lambda)\neq0$.
Note,
furthermore,
that there are notable manifolds
for which vanishing of $\psi(\lambda)$
is equivalent to overtwistedness.
\Def{
	A {\it small Seifert fibred\/} $3$-manifold
	is a Seifert fibred manifold over $S^2$
	with precisely three singular fibres.
}
\Thm{
	(\MatkovicClassificationAlt, Corollary 1.4)
	Let $Y$ be a small Seifert fibred L-space
	and $\lambda$ a contact form on $Y$.
	It follows that $\lambda$ is overtwisted
	if and only if $\psi(\lambda)=0$.
	Furthermore,
	any other tight contact structure $\lambda'$
	is isotopic to $\lambda$
	if only if it be homotopic to $\lambda$.
}
In this situation,
one can reinterpret {\ModPLSpaceLiftContactInvariant}
as reducing a problem in contact topology
to one in homotopy theory.
\Cor{
	Let $G$ be an elementary abelian $p$-group
	for $p$ prime.
	Suppose that $G$ act freely 
	on a small Seifert fibred L-space $Y$
	while preserving some contact form $\lambda$.
	If $\psi(\lambda/G)=0$,
	it follows that $\lambda$ is overtwisted.
	Conversely,
	if $\psi(\lambda/G)\neq0$,
	then
	$\lambda$ is tight
	if and only if
	it be {\it homotopic\/}
	to one of the tight contact structures
	of $Y$.
}
\Rem{
	The tight contact structures on
	small Seifert fibred L-spaces are fully classified
	by
	{\WuLegendrian},
	{\GhigginiLiscaStipsiczClassification},
	{\GhigginiTight}
	and
	{\MatkovicClassification}
}
The author shall finish by noting a certain
constraint on the $d_3$ invariant of a contact structure
which covers a contact structure
with non-vanishing contact invariant.
This result had essentially already appeared in {\RosoSeibergWitten}
except that the condition of a minimal L-space
was unnecessarily assumed there.
\Thm{
	\Label\ConstraintOnDThreeInvOfContactCoveringTight
	If 
	$\dim_{\Field}\HMTilde^*(Y,\SpinCStruct_\lambda;\Field)=1$
	and
	$\psi(\lambda/G)\neq 0$,
	then
	$d_3(\lambda)+{1\over 2}+2h(\SpinCStruct_\lambda)\geq 0$.
}
\Ack{
	The author would like to express his gratitude to S.~Sivek
	for many conversations
	which fundamentally influenced the present \Document.
	He would also like to thank J.~Rasmussen
	for once questioning him
	whether the aforementioned results of {\RosoSeibergWitten}
	could be generalized,
	which was his main motivation for writing this \Document.
	This work was supported by the EPSRC grant EP/W522636/1.
}

\Sect{Borel Floer contact invariant}
Suppose $G$ be a finite group
acting on a contact rational homology $3$-sphere
$(Y,\lambda)$
and fix some compatible $G$-invariant metric $g$ on $Y$.
Denote by $\SpinCStruct_\lambda$
the $\SpinC$ structure
canonically defined by $\lambda$
(vid.\ \PetitSpinCAlt).
This section shall recall the construction
of the equivariant contact invariant
$\psi_G(\lambda)\in\BorelHMTilde_G^*(Y,\SpinCStruct_\lambda)$
and prove a few general results about it
which had not been considered in 
\RosoSeibergWitten.
\par
The definition of the Borel Floer cohomology
is performed via the
{\it Seiberg--Witten Floer spectrum\/}
construction
of \ManolescuSeibergWittenFloer;
therefore,
it shall be necessary to briefly recall that construction as well.
Denote by $\Spinor_\lambda\to Y$
the spinor bundle associated to
the $\SpinC$ structure $\SpinCStruct_\lambda$
and use $\Connex(\det\SpinCStruct_\lambda)$
to denote the affine space
of hermitian connexions on the determinant line bundle
of $\SpinCStruct_\lambda$.
Note that one has an isomorphism
$$
	\Spinor_\lambda
	\cong
	\underline{\C}\oplus(\Ker\lambda)^*
$$
Given this isomorphism,
notice that,
naturally associated to the contact form $\lambda$,
there is an evident nowhere vanishing spinor
$\psi_\lambda\in\Gamma\Spinor_\lambda$
defined to be the constant function $1$.
Likewise,
a connexion $A_\lambda$
can be defined by requiring that the Dirac equation
$$
	\Dirac_{A_\lambda}\psi_\lambda=0
$$
be satisfied.
Define the {\it Seiberg--Witten Coulomb gauge
with respect to the connexion $A_\lambda$\/}
to be the affine space
$$
	\Coulomb(Y,\SpinCStruct_\lambda)
	:=
	\Big(
		A_\lambda
		+
		\Ker(
			\d^*:
			i\Omega^1(Y)
			\to
			i\Omega^0(Y)
		)
	\Big)
	\times
	\Gamma\Spinor_\lambda
	\subset
	\Connex(\det\SpinCStruct_\lambda)
	\times
	\Gamma\Spinor_\lambda.
$$
Notice that,
by regarding $(A_\lambda,0)$
as the zero vector,
one can think of
the affine space $\Coulomb(Y,\SpinCStruct_\lambda)$
as a vector space.
\par
Following {\TaubesWeinsteinII},
the Seiberg--Witten vector field
canonically perturbed by $\lambda$
is defined to be
$$
	\eqalign{
		&\SW_{\lambda,r}:
		\Coulomb(Y,\SpinCStruct_\lambda)
		\to\Tan\Coulomb(Y,\SpinCStruct_\lambda),
		\cr
		&\SW_{\lambda,r}(A,\psi)=
		\left(
			\Hodge{1\over2}\left(
				F_A
				-F_{A_\lambda}
			\right)
			+\Pi(r\tau(\psi))
			-{ir\over2}\lambda
			,\;
			\Dirac_A\psi
			+\xi(\psi)\psi
		\right),
		\cr
	}
$$
where $r>0$ denotes a parameter
which shall be chosen large
according to a few technical requirements
to be explained shortly,
$\Pi$ denotes the $\LTwo$ orthogonal projection
onto $\Ker\d^*$,
$\tau(\psi)$
denotes the inverse via Clifford multiplication
of the spinor bundle endomorphism
$\psi\otimes\psi^*-{1\over2}|\psi|^2$,
and $\xi(\psi):Y\to i\R$ is a complex valued function
defined by requiring that
$\d\xi(\psi)=(1-\Pi)\tau(\psi)$
and $\int_Y\Hodge\xi(\psi)=0$.
\par
Notice that one can split the operator
$\SW_{\lambda,r}$
conveniently
into a linear part
and a compact part as follows.
Write
$\SW_{\lambda,r}=\ell+c$
where
$$
	\ell(A_\lambda+a,\psi)
	:=\left(
		A_\lambda+{1\over2}\Hodge\d a,
		\Dirac_{A_\lambda}\psi
	\right).
$$
It follows that $\ell$
is an elliptic first order linear differential operator.
One can then consider the Sobolev completion
$\Coulomb(Y,\SpinCStruct_\lambda)_k$
of $\Coulomb(Y,\SpinCStruct_\lambda)$
and regard $\ell$ as a Fredholm map
of index zero
$
	\Coulomb(Y,\SpinCStruct_\lambda)_k
	\to
	\Coulomb(Y,\SpinCStruct_\lambda)_{k-1}
$.
Henceforth,
fix $k\geq5$
and bear in mind that the author shall suppress $k$
from the notation.
Given an interval
$I\subset\R$
denote by
$$
	\Coulomb^I(Y,\SpinCStruct_\lambda)
	\subset
	\Coulomb(Y,\SpinCStruct_\lambda)
$$
the span of the eigenvectors of $\ell$
having eigenvalues in the interval $I$.
\par
The construction of the Seiberg--Witten Floer spectrum
can now be made
via use of the theory of \ConleyIsolated.
\Thm{
	(\ManolescuSeibergWittenFloerAlt)
	For sufficiently large $\mu>0$,
	the set
	$$
		S^\mu
		\subset
	\Coulomb^{(\mu,\mu)}(Y,\SpinCStruct_\lambda)
	$$
	consisting of all bounded trajectories of the vector
	field $\SW_{\lambda,r}$
	is an isolated invariant set
	(in the sense of \ConleyIsolatedAlt,
	\S I.1)
	with isolating neighbourhood
	the disk
	$\Disk(\Coulomb^{(\mu,\mu)}(Y,\SpinCStruct_\lambda),R)$
	for some arbitrarily large $R>0$.
}
Given this,
define
the pointed $\U(1)$-space
$\SWF(Y,\SpinCStruct_\lambda,g,\mu)$
to be the $\U(1)$-equivariant Conley index
(vid.\ \ConleyIsolatedAlt, \S III)
of $S^\mu$.
Hence,
one defines
the (metric dependent) Seiberg--Witten Floer spectrum as
$$
	\SWF(Y,\SpinCStruct_\lambda,g)
	:=\Sigma^{-\Coulomb^{(-\mu,0)}(Y,\SpinCStruct_\lambda)}
	\SWF(Y,\SpinCStruct_\lambda,g,\mu)
$$
where $\Sigma^{-V}$
denotes the desuspension functor
by a $\U(1)$-representation $V$
on the equivariant stable homotopy category
$\StabHoSpectra{\U(1)}{(\R^\infty\oplus\C^\infty)}$
(vid.\ \MayEtAlEquivariantAlt, \S XII).
\Thm{
	(\ManolescuSeibergWittenFloerAlt)
	The spectrum
	$\SWF(Y,\SpinCStruct_\lambda,g)$
	is an invariant
	of the triple $(Y,\SpinCStruct_\lambda,g)$
	and,
	moreover,
	only depends on $g$
	up to suspensions.
}
For the present \Document,
it is crucial to note that the spectrum
$\SWF(Y,\SpinCStruct_\lambda,g)$
is really a refinement of the
monopole Floer cohomology
as defined in \KronheimerMrowkaMonopolesAndThreeManifolds.
\Thm{
	(\LidmanManolescuEquivalenceAlt)
	There is an isomorphism
	of $\Q$-graded abelian groups
	$$
		\HMTilde^*(Y,\SpinCStruct_\lambda)
		\cong
		\Homology^{*-n(Y,\SpinCStruct_\lambda,g)}
		(\SWF(Y,\SpinCStruct_\lambda,g)),
	$$
	where
	$n(Y,\SpinCStruct_\lambda,g)$
	is a certain rational number
	defined by Formula (6) of
	\ManolescuSeibergWittenFloer.
}
The monopole contact class
$\psi(\lambda)\in\HMTilde^*(Y,\SpinCStruct_\lambda)$
can be defined in two ways.
On the one hand,
one can argue as
{\KronheimerMrowkaOzsvathSzaboMonopoles}
and define it by counting
Seiberg--Witten monopoles
on $\R\times Y$
with the cylindrical metric on
$(-\infty,-1]\times Y$
and a metric compatible with the symplectization
of the contact form $\lambda$
on $[1,\infty)\times Y$.
On the other hand,
it is possible to degenerate this construction in a sense.
{\TaubesWeinsteinII}
demonstrates that when one perturbs the Seiberg--Witten equations
by the contact form $\lambda$ as shown above,
this count of monopoles on $\R\times Y$
become rather simple as shall be explained next.
\par
After a simple examination of the above formula for
$\SW_{\lambda,r}$,
one finds that the configuration $(A_\lambda,\psi_\lambda)$
is indeed a solution
to the Seiberg--Witten equations $\SW_{\lambda,r}=0$.
The importance of the parameter $r>0$
is that,
if it be chosen large,
the solution $(A_\lambda,\psi_\lambda)$
is non-degenerate
(vid. \TaubesWeinsteinIAlt, Lemma 3.3)
and there are no non-trivial Seiberg--Witten monopoles
on $\R\times Y$
with the cylindrical metric everywhere
for which the positive time limit
be $(A_\lambda,\psi_\lambda)$
(vid.\ \RosoSeibergWittenAlt, Theorem 2.25).
These two properties
imply that the gauge equivalence class of the monopole
$(A_\lambda,\psi_\lambda)$
defines a cocycle in the monopole Floer cochain complex
and,
therefore,
a cohomology class
$\psi(\lambda)\in\HMTilde^*(Y,\SpinCStruct_\lambda)$.
\par
As shown in {\TaubesWeinsteinII}, Proposition 4.3,
the two approaches to the definition of $\psi(\lambda)$
coincide because the count of monopoles on $\R\times Y$
with the cylindrical metric on the negative end
and the symplectization metric on the positive end
consists of a unique gauge equivalence class
whose negative time limit is no other than
$(A_\lambda,\psi_\lambda)$.
However one may choose to define $\psi(\lambda)$,
it turns out to be well defined only up to a sign.
\par
The main result of {\RosoSeibergWitten}
was to refine the definition of $\psi(\lambda)$
to the setting of the spectrum $\SWF(Y,\SpinCStruct_\lambda,g)$
while recovering the cohomology class $\psi(\lambda)$
by a Hurewicz map.
\Thm{
	\Label\CohomotopicalContactInvariant
	(\RosoSeibergWittenAlt, Definition 4.28, Theorem 5.9)
	There exists a $\U(1)$-equivariant cofibre map
	$$
		\Psi(\lambda,g):
		\SWF(Y,\SpinCStruct_\lambda,g)
		\to
		\ContactThom(\lambda,g)
	$$
	where $\ContactThom(\lambda,g)$ is the Thom spectrum
	of a virtual vector bundle over a free $\U(1)$-orbit
	such that
	$\psi(\lambda)$
	corresponds to the class
	$$
		\Psi(\lambda,g)^*\ThomClass_\lambda
		\in\Homology^*_{\U(1)}
		(\SWF(Y,\SpinCStruct_\lambda,g))
		\cong
		\HMTilde^{*+n(Y,\SpinCStruct_\lambda,g)}
		(Y,\SpinCStruct_\lambda)
	$$
	under the isomorphism of \LidmanManolescuEquivalence,
	where $\ThomClass_\lambda$
	is a choice of Thom class for the Thom spectrum
	$\ContactThom(\lambda,g)$.
}
The construction of $\Psi(\lambda,g)$
is performed by reinterpreting
the fact that there are no Seiberg--Witten trajectories
with forward time limit at the configuration
$(A_\lambda,\psi_\lambda)$
in a Conley theoretic manner.
\Prop{
	(\RosoSeibergWittenAlt, Theorem 2.24)
	For sufficiently large
	$\mu>0$ and $r>0$,
	the $\U(1)$-orbit
	of the configuration $(A_\lambda,\psi_\lambda)$
	in $S^\mu$
	is a {\it repeller\/}
	in the sense of \ConleyIsolated,
	Definition 5.1.
}
It is a general construction in Conley theory
that repellers induce {\it connection maps;}
vid.\ \ConleyIsolated, \S III.7.
Hence,
the map
$\Psi(\lambda,g)$
is simply defined
as the desuspension of the connection map
from $\SWF(Y,\SpinCStruct,g,\mu)$
to Conley index
of the $\U(1)$-orbit of $(A_\lambda,\psi_\lambda)$,
which is a Thom spectrum
by non-degeneracy of $(A_\lambda,\psi_\lambda)$.
\par
Perhaps the main advantage of the spectrum apparatus
is the ease of working equivariantly.
Indeed,
it is not difficult to see
that,
in the case at hand,
$G$ acts on the space
$\SWF(Y,\SpinCStruct_\lambda,g,\mu)$
and that its desuspension
$\SWF(Y,\SpinCStruct_\lambda,g)$
becomes a
$(\U(1)\times G)$-spectrum.
When thinking of these
equivariantly,
the author shall denote them by $\SWF_G$.
The key result of
{\LidmanManolescuCovering}
was to establish that,
when $G$ acts {\it freely\/} on $Y$,
the spectra of $Y$ and $Y/G$ are related by
$$
	\Phi^G\SWF_G(Y,\SpinCStruct_\lambda,g)
	=\SWF(Y/G,\SpinCStruct_{\lambda/G},g/G),
$$
where $g/G$ denotes the induced metric on $Y/G$
and $\Phi^G$ denotes the {\it geometric fixed points functor\/}
(vid.\ \MayEtAlEquivariantAlt, \S XVI.3).
The point now is that same relation holds
for the cofibre map $\Psi(\lambda,g)$.
\Thm{
	(\RosoSeibergWittenAlt, \S6)
	The cofibre map $\Psi(\lambda,g)$
	is $G$-equivariant
	for a certain natural $G$-action
	on the Thom spectrum $\ContactThom(\lambda,g)$
	and the induced map on fixed points is
	$\Phi^G\Psi(\lambda,g)=\Psi(\lambda/G,g/G)$.
}
When thinking of $\Psi$ equivariantly,
the author shall denote it by
$$
	\Psi_G(\lambda,g):
	\SWF_G(Y,\SpinCStruct_\lambda)
	\to
	\ContactThom_G(\lambda,g).
$$
Now,
if,
upon passage to cohomology,
one obtains $\psi(\lambda)$
from $\Psi(\lambda,g)$,
one can likewise work in
{\it $G$-equivariant Borel cohomology\/}
to define an equivariant contact class $\psi_G(\lambda)$.
\Def{
	Given a pointed $G$-space $X$,
	define its $G$-equivariant Borel cohomology
	with coefficients in a ring $R$ as
	$$
		\BorelHomology^*_G(X;R)
		:=\Homology^*(X\wedge\E G_+;R),
	$$
	where $\E G_+$ denotes a free contractible $G$-CW complex
	with a disjoint base-point added.
}
\Rem{
	For a $G$-spectrum $X$,
	one can define its Borel cohomology
	by taking a colimit as explained in
	\RosoSeibergWitten, Definition 3.73.
}
\Def{
	Let the {\it $G$-equivariant monopole Borel Floer cohomology\/}
	be
	$$
		\BorelHMTilde_G^*(Y,\SpinCStruct_\lambda;R)
		:=\BorelHomology^{*-n(Y,\SpinCStruct_\lambda,g)}_G
		(\SWF_G(Y,\SpinCStruct_\lambda,g);R).
	$$
}
\Def{
	Given a virtual $G$-vector bundle $E\to B$
	over a $G$-spectrum $B$,
	denote its {\it Thom spectrum\/} by $T(E)$.
	An element $\ThomClass\in\BorelHomology^*_G(T(E);R)$
	is called a {\it Thom class\/}
	if $\BorelHomology^*_G(T(E);R)$
	be a free rank one module over $\BorelHomology^*_G(B;R)$
	with generator $\ThomClass$.
}
\Rem{
	A Thom class need not always exist;
	however,
	if $G$ be an elementary abelian $p$-group
	for $p$ prime
	and the ring of coefficients $R$
	be chosen as the field of $p$ elements,
	then such a Thom class always exists
	by the Thom isomorphism theorem.
}
\Def{
	Assume that a $G$-equivariant Thom class exist
	for the Thom $G$-spectrum $\ContactThom_G(\lambda,g)$
	and denote it by  $\ThomClass_{\lambda,G}$.
	The {\it $G$-equivariant contact element\/} is
	$\psi_G(\lambda):=\Psi_G(\lambda,g)^*\ThomClass_{\lambda,G}$.
}
\Rem{
	In the following sections,
	the above restrictions on $G$ and $R$ shall always be assumed
	so that the Thom class $\ThomClass_{\lambda,G}$ always exist.
}
\Rem{
	Note that one can choose the non-equivariant Thom class
	$
		\ThomClass_{\lambda}
		\in\BorelHomology^*(\ContactThom_G(\lambda,g))
	$
	so that the restriction
	of $\ThomClass_{\lambda,G}$
	to the fibre of the Borel construction
	$
		\ContactThom_G(\lambda,g)
		\to
		\left(\ContactThom_G(\lambda,g)
		\wedge\E G_+\right)/G
	$
	coincide with $\ThomClass_\lambda$.
}
\Rem{
	As with $\psi(\lambda)$,
	except for the case when the ring $R$
	be the field of two elements,
	$\psi_G(\lambda)$ is only defined up to sign
	because the Thom class is only defined up to sign.
}
To finish this section,
the author shall note that
the invariant $\psi(\lambda)$
is homogeneous with respect to the $\Q$-grading
of $\HMTilde^*(Y,\SpinCStruct_\lambda)$
and its grading is given by $d_3(\lambda)+1/2$,
where $d_3$ denotes the three-dimensional
obstruction invariant of hyperplane fields
of \GompfHandlebody,
whose definition shall be recalled next
for the purposes of clarifying the author's conventions.
\Def{
	Let $X$ be a connected
	$4$-manifold-with-boundary
	with an almost complex structure $J$
	such that $(Y,\lambda)$
	be the almost complex boundary of $(X,J)$;
	that is,
	$\Ker\lambda=\Tan Y\cap J\Tan Y$.
	Define
	$$
		d_3(\lambda):=
		{1\over4}\left(
			c_1(X,J)^2-2\chi(X)-3\sigma(X)
		\right).
	$$
	where,
	$c_1(X,J)$ denotes the first Chern class of $\Tan X$
	as a complex vector bundle and,
	for a class $c\in\Homology^2(X;\Z)$,
	the square $c^2\in\Q$
	is defined by taking a preimage
	$\tilde c\in\Homology^2(X,Y;\Q)$
	of $c$
	and setting
	$c^2=(\tilde c\smile\tilde c)[X,Y]$.
}
The $\Q$-grading of the monopole Floer cohomologies
was essentially already noted by
{\FroyshovSeibergWitten}
and is discussed in detail in
\KronheimerMrowkaMonopolesAndThreeManifolds,
\S28.3.
In \LidmanManolescuEquivalence,
Proposition 9.2.1,
it is shown that the grading of the homology
of the spectrum,
$
	\Homology^{*-n(Y,\SpinCStruct_\lambda,g)}
	(\SWF(Y,\SpinCStruct_\lambda,g))
$,
also agrees with the $\Q$-grading of
$\HMTilde^*(Y,\SpinCStruct_\lambda)$.
Hence,
note
(cf.\ \KronheimerMrowkaMonopolesAndThreeManifoldsAlt,
Definition 28.3.1;
\KronheimerMrowkaMonopolesAndContactStructuresAlt,
Theorem 2.4)
that
the contact class $\psi(\lambda)$
has degree
$$
	\deg(\psi(\lambda))
	= d_3(\lambda) + {1\over2}.
$$
It,
therefore,
follows that the $G$-equivariant contact class
has the same degree.

\Sect{Serre Spectral Sequence}
Consider now the behaviour of the contact invariants
in the Serre spectral sequence.
Henceforth,
the author shall always use
the field $\Field:=\Z/p\Z$
as coefficients for the cohomology theories;
this choice shall sometimes be left implicit.
Recall that the second page
of the Serre spectral sequence
is the group cohomology
$\Homology^*(G,\HMTilde^*(Y,\SpinCStruct_\lambda))$,
where $\HMTilde^*(Y,\SpinCStruct_\lambda)$
has a possibly non-trivial $\Z[G]$-module structure
coming from the fact
that it is the (degree-shifted) cohomology
of the $G$-space $\SWF_G(Y,\SpinCStruct_\lambda,g,\mu)$.
In any event,
given any element $\alpha\in\Homology^*(G,\Field)$,
there shall be a class in the second page
$$
	[\alpha,\psi(\lambda)]
	\in\Homology^{\deg x}
	(G,\HMTilde^{\deg\psi(\lambda)}(Y,\SpinCStruct)).
$$
The question is to what this class converges
in the last page.
One might be tempted to say it is $\alpha\psi_G(\lambda)$,
but, because the spectral sequence of a filtered complex
converges to the {\it associated graded module,}
the answer is slightly more involved.
\par
To explain this matter,
introduce a skeletal filtration
$B_k G\subset \B G$,
where $B_k G$ has only cells of dimension up to $k\in\N$.
Then,
the cohomology ring
$\Homology^*(G)\cong\Homology^*(\B G_+)$
is filtered by the ideals
$$
	I_k:=\Ker\Big(
		\Homology^*(\B G_+)
		\to\Homology^*(\B_k G_+)
	\Big)
$$
induced by the inclusion maps.
To simplify the notation,
introduce
$$
	U:=\SWF_G(Y,\SpinCStruct_\lambda,g),
	\quad
	V:=\ContactThom_G(\lambda,g).
$$
Denote by $U_G\to\B G_+$
and $V_G\to \B G_+$ the Borel constructions.
The Borel cohomologies are then filtered by the submodules
$$
	M_k:=I_k\cdot\Homology^*(U_G),
	\quad
	N_k:=I_k\cdot\Homology^*(V_G).
$$
The Serre spectral sequences
of $U_G$ and $V_G$
converge to the associated graded modules
$$
	\AssocGrad\Homology^*(U_G)
	:=\bigoplus_{k\in\N} M_k/M_{k+1},
	\quad
	\AssocGrad\Homology^*(V_G)
	:=
	\bigoplus_{k\in\N} N_k/N_{k+1}.
$$
Consider now the maps
$$
	g_U:
	\Homology^*(\B G_+)\otimes\Homology^*(U_G)
	\to\AssocGrad\Homology^*(U_G),
	\quad
	g_V:
	\Homology^*(\B G_+)\otimes\Homology^*(U_G)
	\to\AssocGrad\Homology^*(U_G),
$$
defined by
$$
	g_U(\alpha\otimes\beta)
	\mapsto
	(\alpha\cdot\beta)M_{\deg\alpha+1}
	\in M_{\deg\alpha}/M_{\deg\alpha + 1}
	\subset
	\AssocGrad\Homology^*(U_G),
$$
and likewise for $g_V$.
\par
Now,
the author can present
the key technical lemma of the present \Document.
\Lem{
	\Label\ConvOfNonEquivInvInSerreSpecSeq
	In the Serre spectral sequence
	of $U$,
	the class
	$[\alpha,\psi(\lambda)]$
	converges to
	$g_U(\alpha\cdot\psi_G(\lambda))$.
}
\Proof{
	Firstly,
	one needs to see that
	$[\alpha,\psi(\lambda)]$
	survives to the last page.
	What matters is naturality
	of the Serre spectral sequence.
	Denote the Serre spectral sequence of
	$U$
	by $E_k$
	and that of $V$
	by $F_k$.
	Next,
	note that $\Psi_G(\lambda,g)$
	induces a map of spectral sequences
	$f_r:F_r\to E_r$;
	such a map intertwines the differentials.
	Because $V$
	is the Thom space
	of a virtual vector bundle over $G$-trivial circle,
	it follows that the $\Field[G]$-module
	$\Homology^*(V)$
	is $G$-trivial.
	Hence,
	on the second page of $F$
	there exists the element
	$$
		\alpha\otimes\ThomClass_\lambda
		\in
		\Homology^{\deg\alpha}(G)
		\otimes
		\Homology^{\deg\psi(\lambda,g)}
		(\ContactThom_G(\lambda,g))
		\cong
		\Homology^{\deg\alpha }
		(G,\Homology^{\deg\psi(\lambda,g)}
		(\ContactThom_G(\lambda,g)))
	$$
	and the spectral sequence $F$
	converges already at page 2.
	It is not difficult to see,
	straight from \CohomotopicalContactInvariant,
	that
	$
		[\alpha,\psi(\lambda)]
		=f_2(\alpha\otimes\ThomClass_\lambda)
	$.
	In particular,
	$[\alpha,\psi(\lambda)]$
	is in the kernel of the differential of the second page
	and any subsequent page.
	This means it survives to $E_\infty$.
	Moreover,
	note that,
	also for any $k\geq 2$,
	$
		[\alpha,\psi(\lambda)]
		=f_k(\alpha,\ThomClass_\lambda)
	$.
	Now,
	because
	$F_2=F_\infty$,
	the spectral sequence $F$
	is completely understood
	and one can check directly that
	$\alpha\otimes\ThomClass_\lambda$
	corresponds to
	$g_V(\alpha\otimes\ThomClass_\lambda)$
	under the isomorphism
	$F_\infty\cong\AssocGrad\Homology^*(V_G)$.
	Hence,
	one sees that
	$[\alpha,\psi(\lambda)]$
	converges to
	$f_\infty(g_V(\alpha\otimes\ThomClass_\lambda))$.
	Furthermore,
	recall that one can choose
	the Thom class
	$\ThomClass_\lambda$
	so that it be the restriction
	to the fibre $V\incl V_G$ of the Borel construction
	of the equivariant Thom class $\ThomClass_{\lambda,G}$.
	This is to say that,
	under the Thom isomorphism
	$
		\Homology^*(U_G)
		\cong
		\Homology^*(G)
		\otimes
		\Homology^*(U)
	$,
	the class $\ThomClass_{\lambda,G}$
	corresponds to
	$1\otimes\ThomClass_\lambda$.
	Therefore,
	one has
	$
		g_V(\alpha\otimes\ThomClass_\lambda)
		=g_V(\alpha\cdot\ThomClass_{\lambda,G})
	$.
	Lastly,
	because the associated graded module construction
	is functorial,
	$
		g_U(\alpha\cdot\psi_G(\lambda))
		=
		f_\infty
		(g_V(\alpha\cdot\ThomClass_{\lambda,G}))
	$.
}
Next,
it is necessary to better understand the map $g_X$.
For that,
recall that the group cohomology of $G$,
if the prime $p$ be odd,
has the following form
$$
	\Homology^*(G)
	=
	\Field[u_1,\cdots,u_n]
	\otimes
	\Exterior_\Field[v_1,\cdots v_n],
$$
where $\deg v_i=1$,
$\deg u_i=2$
and $\Exterior_\Field$
denotes the exterior algebra over the field $\Field$.
On the other hand,
if $p=2$,
one simply has
$$
	\Homology^*(G)
	=
	\Field[u_1,\cdots,u_n]
$$
where $\deg u_i=1$.
One can check that the ideal $I_1\subset\Homology^*(G)$
is the ideal generated by $\{u_i,v_i\mid i=1,\ldots n\}$
if $p\neq 2$
and by $\{u_i\mid i=1,\ldots,n\}$
if $p=2$.
The $k^{\rm th}$ level of the filtration,
$I_k$,
can be expressed as the $k$-fold
product of $I_1$ with itself.
This implies the following.
\Lem{
	\Label\CharOfNonZeroAssocGradEquivInv
	$\psi_G(\lambda)$
	is in the image of the action of one
	of the generators of $\Homology^*(G)$
	in the presentation above
	if and only if
	$g_U(\psi_G(\lambda))=0$.
}
\Proof{
	Note
	that the map
	$g_U$
	restricted to
	$
		\Homology^0(G)
		\otimes
		\BorelHMTilde_G^*(Y,\SpinCStruct_\lambda)
	$
	factors through
	$
		\BorelHMTilde^*_G(Y,\SpinCStruct_\lambda)
		/
		I_1
		\BorelHMTilde^*_G(Y,\SpinCStruct_\lambda)
	$.
	Hence,
	$g_U(\psi_G(\lambda))=0$
	if and only if
	$$
		\psi_G(\lambda)
		\in
		I_1
		\BorelHMTilde^*_G(Y,\SpinCStruct_\lambda).
	$$
	Given the description of the ideal $I_1$ above,
	the result follows.
}
\Cor{
	If $\psi_G(\lambda)$
	is not in the image of one
	of the generators of $\Homology^*(G)$
	in the presentation above,
	then $\psi(\lambda)\neq0$.
}
\Proof{
	Corollary of
	{\ConvOfNonEquivInvInSerreSpecSeq}
	and
	{\CharOfNonZeroAssocGradEquivInv}.
}
It is natural to define
following numerical invariant
which measures how deep $\psi_G(\lambda)$
lies within the filtration
$I_k\cdot\BorelHMTilde_G^*(Y,\SpinCStruct_\lambda)$.
\Def{
	Given the $G$-equivariant contact form $\lambda$
	on the rational homology $3$-sphere $Y$.
	Define the
	{\it filtration invariant\/}
	$\FreeNumericalInvariant_G(\lambda)\in\Z$
	to be the smallest integer
	such that
	$
		\psi_G(\lambda)
		\notin 
		I_{\FreeNumericalInvariant_G(\lambda)+1}
		\BorelHMTilde_G^*(Y,\SpinCStruct_\lambda)
	$.
	In the degenerate case
	$\psi_G(\lambda)=0$,
	set $\FreeNumericalInvariant_G(\lambda):=-1$.
}
\Cor{
	\Label\FreeNumInvVanishesImpliesNonEquivariantNonVanishing
	If $\FreeNumericalInvariant_G(\lambda)=0$
	then $\psi(\lambda)\neq0$.
}

\Sect{Torsion}
This section shall consider,
in some sense,
the opposite type of behaviour of the previous section.
Let $Y$, $\lambda$, $G$ and $\{I_k\}$ be as before.
\Def{
	Suppose that
	$\psi_G(\lambda)\in\BorelHMTilde^*_G(Y,\SpinCStruct_\lambda)$
	be a torsion element
	with respect to the $\Homology^*(G)$-module structure.
	Define the
	{\it torsion invariant\/}
	$\TorNumericalInvariant_G(\lambda)$
	to be the minimum $k\in\N$
	such that $x\cdot\psi_G(\lambda)$
	for all $x$ in the ideal $I_k$.
	On the other hand,
	if $\psi_G(\lambda)$ be a free element,
	define $\TorNumericalInvariant_G(\lambda):=\infty$.
}
This invariant should be seen as a $G$-equivariant analogue of
the invariant of $\sigma$ of \KarakurtContact.
In order to see the r\^ole played by this invariant
in the lifting behaviour of contact structures,
recall firstly the {\it localization theorem\/}
for the Borel Floer cohomology.
\Rem{
	In what follows,
	let $S\subset\Homology^*(G)$
	denote the multiplicative subset
	consisting of the {\it Euler classes\/}
	$e(V)$
	of all $G$-representations $V$
	such that $V^G=0$;
	that is,
	the {\it free\/} $G$-representations.
	For a $\Homology^*(G)$-module $M$,
	denote its localization with respect to $S$
	by $S^{-1}M$.
	Also denote by $S^{-1}$
	the natural map
	$$
		M
		\to
		S^{-1}M
	$$
	Recall that the {\it inclusion of fixed points\/}
	(vid.\ \RosoSeibergWittenAlt, Theorem 6.9;
	cf.\ \LidmanManolescuCoveringAlt),
	induces a map
	$$
		i:
		\BorelHMTilde^*_G(Y,\SpinCStruct_\lambda)
		\to
		\HMTilde^*(Y/G,\SpinCStruct_{\lambda/G})
		\otimes\Homology^*(G).
	$$
	The reader is cautioned that this map
	shifts the gradings by a factor which is
	somewhat difficult to calculate
	in practice
	(cf.\ \RosoSeibergWittenAlt,
	\S8).
}
\Thm{
	\Label\Localization
	(\BaragliaHekmatiEquivariantAlt, Theorem 7.1)
	The map
	$$
		S^{-1}i:
		S^{-1}\BorelHMTilde^*_G(Y,\SpinCStruct_\lambda)
		\to
		\HMTilde^*(Y/G,\SpinCStruct_{\lambda/G})
		\otimes S^{-1}\Homology^*(G),
	$$
	defined by applying the localization functor $S^{-1}$
	to the map $i$,
	is an isomorphism.
}
\Thm{
	\Label\ContactElementNaturalUnderLocalization
	(\RosoSeibergWittenAlt,
	Proposition 7.15)
	The map $S^{-1}i$
	sends $S^{-1}\psi_G(\lambda)$
	to $\psi(\lambda/G)\otimes e(V)$
	for some {\it free\/} virtual representation
	$V$ of $G$;
	in particular,
	the Euler class $e(V)\in\Homology^*(G)$
	does not vanish.
}
\Rem{
	The author should note that
	the result in {\RosoSeibergWitten}
	is stated for the case when $G$ be the cyclic group of order $p$,
	but it is indeed quite immediate
	to generalize it to the case at hand
	of a general elementary abelian $p$-group.
}
\Thm{
	Suppose that $G$ act freely on $Y$.
	If $\TorNumericalInvariant_G(\lambda)=\infty$
	and $\psi(\lambda/G)=0$,
	then $\psi(\lambda)=0$.
}
\Proof{
	By {\Localization}
	and
	{\ContactElementNaturalUnderLocalization},
	if $\psi(\lambda/G)=0$,
	then $S^{-1}\psi_G(\lambda)$
	must also vanish.
	But $S^{-1}\psi_G(\lambda)$
	vanishes if and only if
	there exists $x\in S$
	such that
	$x\cdot\psi_G(\lambda)=0$;
	by examining the $\Homology^*(G)$-module structure
	of $\BorelHMTilde_G^*(Y,\SpinCStruct_\lambda)$,
	one concludes that this happens
	precisely when $\psi_G(\lambda)$ be torsion.
	Hence,
	if $\TorNumericalInvariant_G(\lambda)=\infty$,
	the only option is that $\psi_G(\lambda)=0$.
	To finish off,
	note that,
	by \RosoSeibergWitten, Theorem 7.20,
	there is a map
	$$
		\BorelHMTilde^*_G(Y,\SpinCStruct_\lambda)
		\to
		\HMTilde^*(Y,\SpinCStruct_\lambda)
	$$
	which sends $\psi_G(\lambda)$
	to $\psi(\lambda)$;
	therefore,
	$\psi(\lambda)=0$.
}
If the Borel Floer cohomology
not have any torsion at all,
then,
in particular,
the element $\psi_G(\lambda)$,
if non zero,
must be non-torsion.
\Cor{
	When $\BorelHMTilde_G^*(Y,\SpinCStruct_\lambda)$
	be a free $\Homology^*(G)$-module,
	vanishing of
	$\psi(\lambda/G)$
	implies vanishing of $\psi(\lambda)$.
}
Recall next
the definition of the invariant $h$ of \FroyshovMonopole.
Given a $\SpinC$ structure $\SpinCStruct$ 
on the rational homology $3$-sphere $Y$,
its monopole Floer cohomology
$\HMFrom^*(Y,\SpinCStruct;\Z)$
always takes the form
$\Z[U]\oplus T$
for $T$ a torsion $\Z[U]$-module.
The invariant $h(\SpinCStruct)\in\Q$
is defined so that the grading of the generator
of the rank one free $\Z[U]$-module summand of
$\HMFrom^*(Y,\SpinCStruct;\Z)$
be $-2h(\SpinCStruct)$.
Note,
therefore,
that the simplest possible value
that $\HMTilde^*(Y,\SpinCStruct)$,
now with coefficients in the field $\Field$,
can take is a single copy of $\Field$
concentrated at the degree $-2h(\SpinCStruct_\lambda)$.
\Cor{
	If $\HMTilde^*(Y,\SpinCStruct_\lambda)=\Field$,
	vanishing of
	$\psi(\lambda/G)$
	implies vanishing of $\psi(\lambda)$.
}
An important special case of this condition
is when $\HMTilde^*(Y,\SpinCStruct)=\Field$
for {\it all\/}
$\SpinC$ structures $\SpinCStruct$.
\Cor{
	If $Y$ be a mod $p$ L-space,
	vanishing of
	$\psi(\lambda/G)$
	implies vanishing of $\psi(\lambda)$.
}
Since the author does not know of any example
in which there be $\Homology^*(G)$-torsion in
$\BorelHMTilde_G^*(Y,\SpinCStruct_\lambda)$,
he cannot provide an example in which
$\TorNumericalInvariant_G(\lambda)$ be finite.
\Quest{
	Is it ever the case that
	$\TorNumericalInvariant_G(\lambda)\neq\infty$?
}
This leads to the following more concrete question.
\Quest{
	Do there exist
	$Y$, $G$ and $\lambda$
	for which $\psi(\lambda/G)=0$ but $\psi(\lambda)\neq 0$.
}
The author suspects the answer to both of these question to be positive.

\Sect{Knot Invariants}
Consider a knot $K\subset S^3$
transverse to the standard contact structure
$\xi_{\rm std}:=\Ker\lambda_{\rm std}$.
For a prime $p$,
denote by $\pi:\Sigma_p(K)\to S^3$
the unique connected regular $p$-fold branched cover of $S^3$
with branching locus $K$;
vid. \RolfsenKnots, Chapter 10, \S C,
for a description.
The manifold $\Sigma_p(K)$
is a rational homology sphere;
this fact has been known for a while
(cf.\ \AlexanderBriggsOnTypesAlt)
for a modern proof,
vid.\ \LivingstonSeifert,
Corollary 3.2.
By a well known construction of
\GonzaloBranched,
the manifold
$\Sigma_p(K)$
canonically inherits a contact structure
$\xi(K,p)=\Ker\lambda(K,p)$
from $S^3$ via $\pi$.
Outside a small tubular neighbourhood of $K$,
the contact form $\lambda(K,p)$ is simply given
by the pullback of $\lambda_{\rm std}$
via $\pi$.
Inside the tubular neighbourhood,
one must perform a certain adjustment
so that the contact condition still be satisfied.
Unfortunately,
this adjustment
means that the contact forms $\lambda(K,p)$
and $\lambda_{\rm std}$
are not simply related by the pullback via $\pi$.
Nevertheless,
one can ensure
in the construction of {\GonzaloBranched}
that the contact form $\lambda(K,p)$
be equivariant
with respect to the group $G$ of deck transformations.
This allows one to consider the class
$$
	\psi_G(\lambda(K,p))
	\in\BorelHMTilde^*_G(
		\Sigma_p(K),
		\SpinCStruct_{\lambda(K,p)}
	)
$$
and define numerical invariants of $K$.
\Def{
	Define the {\it filtration invariant\/}
	$\FreeNumericalInvariant(K,p)
	:=\FreeNumericalInvariant_G(\lambda(K,p))$.
}
\Def{
	Define the {\it torsion invariant\/}
	$\TorNumericalInvariant(K,p)
	:=\TorNumericalInvariant_G(\lambda(K,p))$.
}
Currently,
the author has no good methods of computing
either of these invariants
except for the following
general observations.
\Prop{
	Suppose $\Sigma_p(K)$
	be a modulo $p$ L-space,
	then
	$\TorNumericalInvariant(K,p)=\infty$.
	Moreover,
	if $\lambda(K,p)$
	be symplectically fillable,
	then
	$\FreeNumericalInvariant(K,p)=0$.
}
\Cor{
	For a quasi-alternating knot $K$,
	one has
	$\TorNumericalInvariant(K,2)=\infty$.
}
A $G$-equivariant
version of the Heegaard Floer cohomology,
denoted $\HFHat_{G}^*(Y,\SpinCStruct)$
was introduced by \HendricksLipshitzSarkarFlexible.
Although the group $G$ was kept general in principle,
the authors focus on the case $G=\Z/2\Z$
and use coefficients in the field $k=\Z/2\Z$.
It was first conjectured by
{\BaragliaHekmatiEquivariant}
that there might be an isomorphism between
$\HFHat_{G}^*(Y,\SpinCStruct)$
and
its Seiberg--Witten counterpart
$\BorelHMTilde_G^*(Y,\SpinCStruct_\lambda)$.
In \KangZModTwoEquivariant,
an invariant of a transverse knot $K$,
$$
	c_{\Z/2\Z}(\xi(K,2))\in
	\HFHat^*_{\Z/2\Z}
	(\Sigma_2(K),\SpinCStruct_{\lambda(K,2)})
$$
is defined by adapting the construction
of the Heegaard Floer contact element
$c(\xi(K,2))\in\HFHat^*(\Sigma_2(K),\SpinCStruct_{\lambda(K,2)})$
to the equivariant setting.
Given the parallels between the Heegaard Floer
and Seiberg--Witten Floer theories,
it is natural to conjecture the following.
\Conj{
	Under the {\it conjectural\/}
	isomorphism
	$
		\HFHat_{G}^*(Y,\SpinCStruct)
		\cong
		\BorelHMTilde_G^*(Y,\SpinCStruct_\lambda)
	$,
	the class $c_{\Z/2\Z}(\xi(K,2))$
	corresponds to $\psi_{\Z/2\Z}(\lambda(K,2))$.
}

\Sect{Lifting contact classes}
This section is dedicated to finishing the proof of the main results.
The author starts by recalling a result
which follows from \LidmanManolescuCovering;
a proof for this result shall be presented here
using the Borel Floer cohomology.
In \LidmanManolescuCovering,
this result is achieved by invoking
the inequality of \SmithTransformations.
It is well known that
the Smith theory can be recovered using
Borel cohomology;
vid.\ \BorelNouvelle.
The author shall use the idea of that proof
in the particular case at hand
since it shall prove useful
to keep it in mind later.
\par
Let $Y$, $G$ and $\lambda$ be as in the previous sections
and assume that $G$ act freely on $Y$.
The author reminds the reader
that he is implicitly using Floer cohomologies
and group cohomologies
with coefficients in the field $\Field$ of $p$ elements.
\Lem{
	\Label\BorelFloerCohomologyForRankOneHMTilde
	If $\dim_\Field\HMTilde^*(Y,\SpinCStruct_\lambda)=1$,
	$
		\BorelHMTilde_G^*(Y,\SpinCStruct_\lambda)
		=
		\Homology^{*+2h(\SpinCStruct_\lambda)}(G)
	$.
}
\Proof{
	Start with the Serre spectral sequence
	for $\BorelHMTilde^*_G(Y,\SpinCStruct_\lambda)$.
	The second page
	looks like the following.
	$$
	\vrule
	\vbox{%
		\hbox{%
			\hskip0.5em
			\vbox{\halign{%
				\hfil$#$\hfil&&\quad\hfil$#$\hfil\cr
				\vdots & & \vdots & \vdots & \vdots & \cr
				0 & \cdots & 0 & \Field & 0 & \cdots\cr
				0 & \cdots & 0 & \Field & 0 & \cdots\cr
				0 & \cdots & 0 & \Field & 0 & \cdots\cr
				0 & \cdots & 0 & \Field & 0 & \cdots\cr
			}}
			\hskip0.5em
		}
		\vskip 0.5em
		\hrule
	}
	$$
	Hence,
	the differentials already all vanish
	for degree reasons.
	Since the coefficients of the cohomologies
	are in a field,
	the result follows. 
}
\Cor{
	(cf.\ \LidmanManolescuCoveringAlt, Corollary 1.2)
	If $\dim_\Field\HMTilde^*(Y,\SpinCStruct_\lambda)=1$,
	then $\dim_\Field\HMTilde^*(Y/G,\SpinCStruct_{\lambda/G})=1$.
}
\Proof{
	By {\BorelFloerCohomologyForRankOneHMTilde}
	and {\Localization} (localization),
	it follows that
	$$
		S^{-1}\BorelHMTilde_G^*(Y,\SpinCStruct_\lambda)
		\cong
		S^{-1}\Homology^*(G)
		\cong
		\HMTilde^*(Y/G,\SpinCStruct_\lambda)
		\otimes_\Field
		S^{-1}\Homology^*(G),
	$$
	where the isomorphisms here
	do not preserve the gradings.
	In turn,
	this implies that
	$\HMTilde^*(Y/G,\SpinCStruct_\lambda)$
	must be of dimension precisely $1$.
}
The next result shows that the numerical invariant
$\FreeNumericalInvariant_G$
is easily computable in terms of the $\Q$-grading
of $\HMTilde^*(Y,\SpinCStruct_\lambda)$
in the simplest case.
\Thm{
	\Label\ComputedFreeNumericalInvariantSpinCLSpace
	If $\dim_\Field\HMTilde^*(Y,\SpinCStruct_\lambda)=1$
	and $\psi(\lambda/G)\neq 0$,
	then
	$
		\FreeNumericalInvariant_G(\lambda)
		=d_3(\lambda)+1/2
		+2h(\SpinCStruct_\lambda)
	$.
}
\Proof{
	The ideal $I_k\subset\Homology^*(G)$
	can be defined as the ideal consisting of
	elements with degree no less than $k$.
	{\BorelFloerCohomologyForRankOneHMTilde}
	therefore implies that the filtration
	$I_k\cdot\BorelHMTilde^*_G(Y,\SpinCStruct_\lambda)$
	also consists of elements
	of degree bounded below by
	$k$ minus the minimal degree
	for which $\BorelHMTilde^*_G(Y,\SpinCStruct_\lambda)$
	not vanish;
	that is,
	$-2h(\SpinCStruct_\lambda)$.
	Hence,
	it follows that
	$
		\FreeNumericalInvariant_G(\lambda)
		=
		\deg\psi_G(\lambda)
		+2h(\SpinCStruct_\lambda)
	$.
}
\Eg{
	The contact topology of the lens spaces $L(n,m)$
	is well understood
	according to
	{\HondaClassification}
	or
	{\GirouxStructures}.
	Given $n>m>0$
	such that $\gcd(n,m)=1$,
	write its continued fraction expansion as
	$$
		-{n\over m}=
		r_0
		-{1\over{
			r_1-{
				1\over{
					r_2-\cdots
					{1\over r_k}
				}
			}
		}}.
	$$
	Then,
	there is a total of
	$\Pi_{i=0}^k|r_i+1|$
	isotopy classes of tight contact structures on
	$L(n,m)$.
	Now,
	consider some prime $p$
	and suppose that it divide $n$.
	Then,
	there is a regular cyclic
	covering $\pi:L(n/p,m)\to L(n,m)$
	of degree $p$.
	One can then pullback via $\pi$
	a tight contact form $\lambda$ on $L(n,m)$
	to obtain a $G$-equivariant
	contact form $\pi^*\lambda$ on $L(n/p,m)$.
	Note that
	$\psi(\lambda)\neq0$
	because,
	according to
	{\HondaClassification},
	all tight contact structures on lens spaces
	are Stein fillable
	and,
	therefore,
	their contact classes are non-vanishing
	by {\OzsvathSzaboContact}.
	Since there is a map sending
	$\psi_G(\pi^*\lambda)$
	to $\psi(\lambda)$,
	it must be the case that
	$\psi_G(\pi^*\lambda)\neq0$ as well.
	Furthermore,
	lens spaces are L-spaces,
	and,
	therefore,
	modulo $p$ L-spaces;
	hence,
	one can apply
	{\ComputedFreeNumericalInvariantSpinCLSpace}
	to compute the filtration invariant.
	$$
		\FreeNumericalInvariant_G(\pi^*\lambda)
		=
		d_3(\pi^*\lambda)+{1\over2}
		+2h(\SpinCStruct_{\pi^*\lambda}),
	$$
	Note that,
	in this family of examples,
	one obtains a variety of cases in which
	$\FreeNumericalInvariant_G(\pi^*\lambda)>0$.
	Indeed,
	according to {\HondaClassificationII},
	as long as there be at least two distinct integers
	$k_0,k_1\in\{0,\ldots,k\}$
	such that $r_{k_0}>-2$
	and $r_{k_1}>-2$,
	there exists a value of $p$
	for which its corresponding lift
	to $L(m/p,n)$ shall be overtwisted
	for all but two of the tight contact structures
	of $L(m,n)$.
	Suppose $m$, $n$ and $p$
	be as such.
	Then,
	$\psi(\pi^*\lambda)=0$
	because $\pi^*\lambda$ overtwisted.
	If it were the case that
	$\FreeNumericalInvariant_G(\pi^*\lambda)=0$,
	then $\psi(\pi^*\lambda)$ would be non-vanishing
	by
	\FreeNumInvVanishesImpliesNonEquivariantNonVanishing.
}
\Thm{
	If $\dim_\Field\HMTilde^*(Y,\SpinCStruct_\lambda)=1$
	and $\psi(\lambda/G)\neq 0$,
	then $\psi(\lambda)\neq 0$
	if and only if
	$d_3(\lambda)+{1\over2}+2h(\SpinCStruct_\lambda)=0$.
}
\Proof{
	Follows from combining
	{\FreeNumInvVanishesImpliesNonEquivariantNonVanishing}
	with
	{\ComputedFreeNumericalInvariantSpinCLSpace}.
}
\Cor{
	If $Y$ be a mod $p$ L-space
	and $\psi(\lambda/G)\neq 0$,
	then $\psi(\lambda)\neq 0$
	if and only if
	$d_3(\lambda)+{1\over2}+2h(\SpinCStruct_\lambda)=0$.
}
The author now concludes
by providing the proof,
essentially unchanged from \RosoSeibergWitten,
of \ConstraintOnDThreeInvOfContactCoveringTight.
\Thm{
	If 
	$\HMTilde^*(Y,\SpinCStruct_\lambda)=1$
	and
	$\psi(\lambda/G)\neq 0$,
	then
	$d_3(\lambda)+{1\over 2}+2h(\SpinCStruct_\lambda)\geq 0$.
}
\Proof{
	If $\psi(\lambda/G)\neq 0$,
	then $\psi_G(\lambda)\neq 0$
	by {\Localization} (localization).
	Now,
	recall that
	$\deg(\psi_G(\lambda))=d_3(\lambda)+{1\over2}$.
	Meanwhile,
	by {\BorelFloerCohomologyForRankOneHMTilde},
	$\BorelHMTilde_G^q(Y,\SpinCStruct_\lambda)=0$
	if $q<-2h(\SpinCStruct_\lambda)$.
}
According to {\EliashbergClassification},
the overtwisted contact structures
on $S^3$
are classified by their $d_3$ invariants.
On $S^3$,
the $d_3$ invariant takes values in $\Z+{1\over2}$.
Hence,
one has the following immediate corollary.
\Cor{
	None of the overtwisted contact structures
	on $S^3$
	whose $d_3$ invariants are strictly lower
	than $-1/2$
	can cover a tight lens space.
}

\vfill
\eject
\Sect{Bibliography}
\Bib
\bye